%% file: main.tex
\documentclass[review, 3p, onecolumn]{elsarticle}
\usepackage[expansion=false]{microtype}
\usepackage{graphicx}
\usepackage{mathrsfs}
\usepackage{amssymb}
\usepackage{amsmath}
\usepackage[T1]{fontenc}
\usepackage{booktabs}
\usepackage{longtable}
\usepackage{siunitx} 
\usepackage{setspace}
\usepackage{amsthm} 
\usepackage{lineno}
\usepackage{listings}
\usepackage{multirow}
\usepackage{lscape}
\usepackage{ bbold }
\usepackage{ esint }
\usepackage{algorithm}
\usepackage{algpseudocode}
\usepackage{enumitem}

\usepackage{subcaption}  

\newcommand{\Ti}{T_{\mathrm{i}}}         
\newcommand{\To}{T_{\mathrm{o}}}         
\newcommand{\Ts}{T_{\mathrm{s}}}           
\newcommand{\Tout}{T_{\mathrm{out}}}
\newcommand{\Tin}{T_{\mathrm{in}}}
\newcommand{\Toutsys}{T_{\text{out,sys}}}

\newcommand{\Cf}{C_{\mathrm{f}}}           
\newcommand{\mdot}{\dot{m}}                
\newcommand{\cp}{c_{p}}                    
\newcommand{\qbhe}{q_{\text{bhe}} }
\newcommand{\qbheL}{q_{\text{bhe}, L} }
\newcommand{\qsys}{q_{\text{sys}}}
\newcommand{\Nb}{N_{b} }
\newcommand{\Deltaneigh}{\Delta_{\text{neigh}}}
\newcommand{\Deltabhe}{\Delta_{\text{bhe}}}

\newcommand{\Rdelta}[1]{R_{#1}^{\Delta}}

\newcommand{\Zmid}{\mathcal{Z}_{\text{mid}}}

\newcommand{\pder}[2]{\frac{\partial #1}{\partial #2}}

\usepackage[colorlinks=True]{hyperref}
\usepackage{xurl} 
\biboptions{sort&compress}

\usepackage{array}
\newcolumntype{L}[1]{>{\raggedright\let\newline\\\arraybackslash\hspace{0pt}}m{#1}}
\newcolumntype{C}[1]{>{\centering\let\newline\\\arraybackslash\hspace{0pt}}m{#1}}
\newcolumntype{R}[1]{>{\raggedleft\let\newline\\\arraybackslash\hspace{0pt}}m{#1}}

\def\equationautorefname~#1\null{(#1)\null}

\usepackage{framed} 
\usepackage{multicol} 
\usepackage{nomencl} 
\makenomenclature
\renewcommand*\nompreamble{\begin{multicols}{2}}
\renewcommand*\nompostamble{\end{multicols}}
\newtheorem{proposition}{Proposition}
\newtheorem{remark}{Remark}
\newtheorem{lemma}{Lemma}

\newtheorem{theorem}{Theorem}
\newtheorem{definition}{Definition}

\usepackage{etoolbox}
\renewcommand\nomgroup[1]{%
  \item[\bfseries
  \ifstrequal{#1}{L}{Latin letters}{%
  \ifstrequal{#1}{M}{Greek letters}{%
  \ifstrequal{#1}{F}{Calligraphic letters}{%
  \ifstrequal{#1}{S}{Subscripts and superscripts}{}}}}%
]}

\makeatletter
\let\els@orig@ps@pprintTitle\ps@pprintTitle
\def\ps@pprintTitle{\els@orig@ps@pprintTitle
  \let\els@orig@oddfoot\@oddfoot
  \def\@oddfoot{\hbox to\textwidth{\els@orig@oddfoot\hss}}%
  \let\@evenfoot\@oddfoot}
\makeatother

\begin{document}

\begin{frontmatter}

\title{Fast Gradient-Based Framework for Spatial and Length Optimization of Geothermal Borehole Fields}

\pdfstringdefDisableCommands{%
  \def\corref#1{}%
}

\author[MO]{Oliver Heinzel}
\author[ENS]{Smajil Halilovic\corref{correspondingauthor}}
\ead{smajil.halilovic@tum.de}
\cortext[correspondingauthor]{Corresponding author}
\author[ENS]{Thomas Hamacher}
\author[MO]{Michael Ulbrich}

\address[MO]{Technical University of Munich, Chair of Mathematical Optimization, Germany}
\address[ENS]{Technical University of Munich, Chair of Renewable and Sustainable Energy Systems, Germany}

\begin{abstract}
Common analytical design tools for Borehole Heat Exchanger (BHE) fields size the total drilling length for a given layout but do not optimize the positions of the boreholes. We introduce a framework (OptiGeoS) for the simultaneous, gradient-based optimization of the borehole length and the continuous borehole positions under groundwater advection. It couples the Moving Finite Line Source model with local BHE equations and evaluates the temporal convolutions by Fast Fourier Transforms. The gradients with respect to the length and the positions are computed analytically. A flow-dependent Centroidal Voronoi Tessellation confines each borehole to its own cell, which enforces the spatial separation and provides advection-informed layouts.

Under purely conductive conditions, the fluid temperature extremes of the forward model for a single borehole lie between those of EED and GHEtool, which differ by $0.32\,\mathrm{K}$. For a field of 120 boreholes, the hourly outlet temperature agrees with that of EED with an RMSE of $0.28\,\mathrm{K}$, while the peak temperatures deviate by up to $2\,\mathrm{K}$. For 25 boreholes simulated hourly over 20 years, one forward simulation takes less than $0.5\,\mathrm{s}$ and a full optimization about $8\,\mathrm{min}$ on a laptop CPU.
In a case study with groundwater advection, the optimized layout requires $16.2\,\%$ less borehole length than a fixed grid. It reduces the length increase caused by the advected thermal plume from $20.9\,\%$ to $1.4\,\%$.
\end{abstract}

\begin{keyword}
Borehole Heat Exchanger \sep
Moving Finite Line Source \sep
Gradient-Based Optimization \sep
Fast Fourier Transform \sep
Groundwater Advection\sep
Centroidal Voronoi Tessellation
\end{keyword}

\end{frontmatter}


\section{Introduction}

\input{sections/introduction}

\section{Coupled Thermal Modeling}
We use a two-scale approach as in \cite[Sec.~2]{RUIZCALVO2016720}. The local model treats the borehole as quasi-steady, following \cite[Sec.~5.1]{CLAESSON1988509} with the extensions of \cite[Sec.~2]{Diersch_resistance} and \cite[Secs.~2.1 and~2.4]{CIMMINO2016137}. The global model describes heat transport in the surrounding subsurface \cite{MFLS}.
\subsection{Local BHE Physics} \label{subsec: local bhe}
\input{sections/local_bhe_physics}

\subsection{Global Soil Physics} \label{subsec: global soil}

\input{sections/MFLS.tex}

\subsection{Temporal Coupling}

\input{sections/temporal_coupling}

\subsection{Lumping and Interaction Split} \label{sec: lumping split}

\input{sections/spectral_decomposition}

\subsection{Numerical Implementation via Fast Fourier  Transform}
    
\input{sections/fft}

\section{Gradient-Based Optimization}

\subsection{Constraints} \label{sec: Constraints}
    \input{sections/constraints}

\subsection{Optimization Problem}
\input{sections/Problem_definitions}


\subsection{Spatial Distribution via Centroidal Voronoi Tessellations} \label{section: spatial Distribution via CVT}

\input{sections/Spatial_strategy}

\section{Benchmarking and Verification} \label{sec: validation}


\input{sections/forward_model_validation}

\section{Case Study} \label{sec: case study}

\input{sections/case_study}

\section{Conclusion and Limitations}

\input{sections/conclusion}

\section*{Data and code availability}
The code, input data and scripts that reproduce all results of this paper are available in the \href{https://github.com/oli-heinzel/Fast-Gradient-Based-Framework-for-Spatial-and-Length-Optimization-of-Geothermal-Borehole-Fields}{GitHub repository} (MIT licence).

\section*{Declaration of competing interest}
The authors declare that they have no known competing financial interests or personal relationships that could have appeared to influence the work reported in this paper.

\section*{Declaration of generative AI and AI-assisted technologies in the writing process}
During the preparation of this work, the authors used Claude (Anthropic; model Claude Sonnet~5.5) in the following ways. (i)~Mathematical content: Claude was asked to draft a proof of Theorem~\ref{thm: existence} (Appendix~\ref{appendix: proof thm1}), which the authors checked step by step and revised. For double-checking the authors' analytical derivations, Claude was used to generate Mathematica notebooks that the authors applied to the proofs and closed-form expressions. (ii)~Software: Claude assisted with writing and debugging the Python code that produces all results. Claude set up the case studies, i.e.\ the configuration files, input data and run scripts. (iii)~Manuscript text: Claude was used for language editing, including rewording the authors' proofs for clarity, for checking the consistency of notation, cross-references and reported numbers, and for reviewing drafts and pointing out errors and gaps. The scientific content, the models, the problem formulation and the numerical experiments are the authors' own. After using this tool, the authors reviewed and edited the content as needed and take full responsibility for the content of the published article.

\section*{Acknowledgement}
\label{sec:Acknowledgement}
The work presented in this paper was supported by the German Federal Ministry for Economic Affairs and Climate Action (BMWK) within the scope of the research project OptiGeoS (01256700/1).

\bibliography{bib}

\appendix

\section{Existence, Uniqueness and Representation of the Temperature Change} \label{appendix: proof thm1}

\input{sections/Appendix_proof_mfls}

\section{Closed Form of the Integrated Vertical Response} \label{appendix: closed V}

\input{sections/Appendix_closed_form_V}

\section{Closed Form of Voronoi Centroids} \label{appendix: closed centroid}

\input{sections/Appendix_closed_form_centroid}

\section{Analytical Gradients}
\label{Appendix_Jacobians}

\input{sections/Appendix_derivatives}



\input{sections/nomenclature}

\end{document}

%% file: sections/introduction.tex
Heat is the largest energy end use and accounts for almost half of global final energy consumption \cite{IEA2021Renewables}. Its decarbonization is part of the transition toward carbon neutrality \cite{IEA2021NetZero}.
Ground Source Heat Pumps (GSHPs) are an efficient, low-carbon alternative to fossil-fuel-based heating and cooling systems \cite{BAYER20121256, YANG201016}.
Using the stable temperature of the shallow subsurface, GSHPs can provide base-load heating and cooling for individual buildings and for district heating and cooling networks \cite{BUFFA2019504}. Borehole Heat Exchangers (BHEs) exchange the heat between the heat pump and the ground.
The deployment of BHE fields is limited by their high initial capital cost. The ground heat exchangers account for more than $30\,\%$ of the total cost of a GSHP system and are the largest contributor to its cost premium over conventional systems \cite{Liu2019}.
Over-sizing the total drilling length by a few percent therefore increases the cost of the ground heat exchangers by about the same fraction. Design tools should minimize the total drilling length while keeping the subsurface thermally sustainable.

Analytical models for BHE fields are reviewed in \cite{LI2015178, KUNTUAROVA2024132189}. \citet[Sec.~4]{CIMMINO2013401} developed a method based on the Finite Line Source that computes thermal response functions for boreholes that are not regularly spaced. This removed the dependence on pre-calculated libraries and allowed the analysis of non-grid layouts. Building on this, \citet[Sec.~2.2]{ROSE2026128565} introduced a framework that constructs both short- and long-term g-functions \cite{CIMMINO2014641}. It is tailored to the thermal interactions of parallel-connected BHE systems.
Common design tools use only part of this flexibility. Sizing tools such as Earth Energy Designer (EED) \cite{etde_573620} and GHEtool \cite{ghe_tool} determine the required borehole length for a given layout. GHEtool also accepts custom borehole positions, but neither tool optimizes the positions together with the length. To our knowledge, neither tool accounts for groundwater advection, although advection shifts the thermal plumes downstream and changes the long-term performance \cite[Sec.~4.2]{TYEGINGRAS2014354}.

\citet{EGIDI2023114957} optimized the BHE positions with an analytical heat conduction model to reduce the interference between neighboring BHEs, without groundwater flow. \citet{Beck2013} optimized the positions together with the loads of the individual boreholes in conduction-dominated ground by evolutionary algorithms and linear programming. In their study, optimizing the positions alone reduced the ground temperature change by about $12\,\%$, almost as much as optimizing the loads. De~Paly et~al.~\cite{dePaly2012} optimized the individual loads by linear programming, and \citet{HechtMendez2013} extended this approach to groundwater flow using moving line sources. In their field of $25$ boreholes, the optimized load patterns align with the flow direction, and the benefit of the optimization decreases with increasing flow velocity. None of these studies optimizes the borehole length. Advection can be included analytically by the Moving Finite Line Source (MFLS) \cite{MFLS, TYEGINGRAS2014354}, which we use here.
Numerical simulators based on finite element methods, such as FEFLOW \cite{Diersch.2014} or OpenGeoSys \cite{ogs:6.5.6}, can model heterogeneous domains with high physical fidelity. However, a single simulation of a BHE field over several decades is computationally expensive, and an optimization requires many such simulations. These suites also often lack automated optimization routines.

We propose a framework (OptiGeoS) for the simultaneous optimization of the BHE positions and the borehole length. It couples an analytical forward model that includes groundwater advection with a gradient-based optimizer that handles the thermal and regulatory constraints. The boreholes can thus be placed with respect to the flow direction. Unlike metaheuristics such as genetic algorithms \cite{genetic} or the evolutionary algorithms of \citet{Beck2013}, it uses analytical gradients of the model. \citet{cimmino2014effects} restricted the optimization to a grid, and \citet{BAYER2014445} and \citet{YANG2026141111} used discrete mixed-integer formulations for placement and scheduling. Our framework instead treats the borehole coordinates as continuous variables. It builds on our earlier work \cite{inproceed_1stpaper}, which placed the boreholes first by a modified Lloyd algorithm that maximizes the spacing, and then optimized only the common length $L$ for these fixed positions, with a conduction-only finite line source model. Here, the positions and $L$ are optimized simultaneously with analytical gradients, and groundwater advection is included.

The geometric constraints are handled with Centroidal Voronoi Tessellations (CVT) \cite{spatial_tessellations}, \cite[Secs.~2.1, 2.5 and~2.6]{CVT_application}. The CVT partitions the convex property into cells that are shrunk by half the minimum borehole spacing, and each BHE moves continuously within its own cell. This is more restrictive than pairwise distance constraints, but it replaces these non-convex constraints by linear inequalities.
We solve the optimization problem with the trust-region method \texttt{trust-constr} of SciPy \cite{2020SciPy-NMeth}, which returns locally optimal layouts and borehole lengths.

%% file: sections/local_bhe_physics.tex
For a BHE of length $L$, the coupled heat balance on the inlet $\Ti(z,t)$ and outlet $\To(z,t)$ fluid temperatures reads \cite[Sec.~8.1, Eq.~(104)]{Diersch_pt1}
\begin{align}
    -\Cf \pder{\Ti}{z} &= \frac{\Ti - \Ts}{\Rdelta{1}} + \frac{\Ti - \To}{\Rdelta{12}} \label{eq:balance_in} \\
    \Cf \pder{\To}{z} &= \frac{\To - \Ts}{\Rdelta{2}} + \frac{\To - \Ti}{\Rdelta{12}}, \label{eq:balance_out}
\end{align}
where $z\in[0,L]$ is the depth below the surface and $\Ts(z,t)$ is the borehole wall temperature. The heat capacity rate of one borehole is $\Cf = \mdot \cp$ $[\unit{\watt\per\kelvin}]$, with mass flow rate $\dot{m}$ and specific heat capacity $c_p$. For $\Nb$ boreholes connected in parallel with the total mass flow rate $\dot{m}_{\mathrm{tot}}$ and a uniform flow distribution \cite[Sec.~1]{ROSE2026128565}, $\Cf = \dot{m}_{\mathrm{tot}}\,c_p/\Nb$.

The thermal resistances $\Rdelta 1$, $\Rdelta 2$ between the pipes and the borehole wall and $\Rdelta{12}$ between the pipes follow from the resistance network of \citet[Sec.~2]{Diersch_resistance}; see \cite[Secs.~5 and~8.3]{Diersch_pt1} for different BHE geometries (1U, 2U, CXA, CXC). For the 1U and 2U geometries, we compute the borehole resistance $R_b$ and the internal resistance $R_a$ with the multipole method \cite{CLAESSON2011} instead of the grout correlations of \citet[Sec.~2]{Diersch_resistance}. Both resistances include the convective resistance in the pipes \cite[Eq.~(32)]{Diersch_pt1}. The network is calibrated to reproduce $R_b$ and $R_a$ exactly, and $\Rdelta 1$, $\Rdelta 2$, $\Rdelta{12}$ follow from it \cite[Sec.~8.3]{Diersch_pt1}.
With the bottom boundary condition $\Ti(L,t) = \To(L,t)$, the outlet temperature $\Tout(t) = \To(0,t)$ follows via the Laplace transform \cite[Sec.~8.2]{Diersch_pt1}.

We assume that each BHE has the same extraction rate $\qbhe = \qsys/\Nb$ $[\unit{\watt}]$. The energy balance then gives $\qbhe(t) = \Cf (\Tout(t) - \Tin(t))$ \cite[Sec.~6.4]{hellstroem_thesis}. Together with the uniform heat flux along each borehole in the soil model (Section~\ref{subsec: global soil}), this is the boundary condition BC-I of \citet[Sec.~1]{ROSE2026128565}. With a common inlet temperature and the same $\Cf$, this boundary condition would strictly imply identical outlet temperatures. In our model, the wall temperatures $\Ts$ differ between the boreholes, so \eqref{eq:outlet_final} yields different outlet temperatures $T_{\text{out}}^k$, which are averaged in \eqref{eq:outlet_mean}. Under groundwater advection, upstream and downstream boreholes would in reality also differ in their extraction rates. Since the optimizer exploits these differences in the wall temperatures, the optimized layouts should be read with this approximation in mind.

Substituting the heat load and solving for $\Tout$ yields \cite[Sec.~8.3]{Diersch_pt1} (see also \cite[App.~M.4.3]{Diersch.2014})
\begin{equation} \label{eq:outlet_final}
\Tout(t) = \psi(L) \int_0^L \Ts(\xi, t)\, w_L(\xi)\, d\xi - \frac{\qbhe(t) (\psi(L)  -1 )}{\Cf},
\end{equation}
with the weight $w_L(\xi) := N_W(L-\xi)/D_W(L)$ and, in the notation of \cite[Eq.~(109)]{Diersch_pt1},
\begin{equation} \label{sys: constants}
\begin{aligned}
\beta_j &= \frac{1}{\Rdelta j \Cf},\quad \beta = \frac{\beta_2 - \beta_1}{2}, \quad \gamma = \sqrt{\frac{(\beta_1 + \beta_2)^2}{4} + \beta_{12}(\beta_1 + \beta_2)}, \quad \nu = \frac{\beta_1 + \beta_2}{2 \gamma}, \quad M_{1,2} = \beta \pm \gamma,\\
\psi(z) &= \frac{1}{2 \nu} \coth(\gamma z) + \frac{1}{2},\\
    N_W(z) &= C_1e^{M_1 z} +C_2 e^{M_2 z}, \qquad C_1 = \nu M_1, \quad C_2 = -\nu M_2,\\
    D_W(z) & = D_1 e^{M_1 z} + D_2 e^{M_2 z}, \qquad D_1 = \frac{1 + \nu}{2}, \quad D_2 = \frac{1 - \nu}{2}.
\end{aligned}
\end{equation}

%% file: sections/MFLS.tex
We consider a constant velocity field $\mathbf{v} = (v_x, v_y, 0)^\top \in \mathbb{R}^3$ that represents steady-state groundwater flow. The temperature change $u(\mathbf{x}, t) := T(\mathbf{x}, t) - T_0(z)$ with respect to the undisturbed profile
\begin{equation} \label{eq: T0 linear}
T_0(z) = T_{\text{surface}} + g_{\text{geo}}\, z
\end{equation}
with a constant geothermal gradient $g_{\text{geo}}$ is governed by the linear advection-diffusion equation on $\mathbb{R}^3$
\begin{equation} \label{eq: advection diffusion}
\begin{aligned}
    \partial_t u + \mathbf{v} \cdot \nabla u - \alpha \Delta u &= \frac{\alpha}{\lambda} f(\mathbf{x}, t), && \mathbf{x} \in \mathbb{R}^3, \, t > 0,\\
    u(\mathbf{x}, 0) &= 0, && \mathbf{x} \in \mathbb{R}^3.
\end{aligned}
\end{equation}
Here, $\alpha \ [\unit{\meter\squared\per\second}]$ denotes the thermal diffusivity and $\lambda \ [\unit{\watt \per \meter \per \kelvin}]$ the thermal conductivity of the soil. The velocity $\mathbf{v}=\frac{(\rho c)_w}{(\rho c)_m}\mathbf{v}_{\mathrm D}$ is the effective thermal velocity \cite{MFLS}, i.e.\ the Darcy flux $\mathbf{v}_{\mathrm D}$ scaled by the ratio of the volumetric heat capacities of the water and of the saturated soil, with $(\rho c)_m=\lambda/\alpha$. The source term $f$ represents the $\Nb$ BHEs as finite line sources with the extraction rate $\qbhe$ of Section~\ref{subsec: local bhe}, distributed uniformly along the borehole,
\begin{equation} \label{eq: source f}
    f(\mathbf{x}, t) = -\sum_{k=1}^{\Nb} \qbheL(t) \delta(x-p^k_x)\delta(y-p^k_y) \left( \mathbb{1}_{[0, L]}(z) - \mathbb{1}_{[-L, 0]}(z) \right),
\end{equation}
where the minus sign makes $\qbhe>0$ denote heat extraction, $\qbheL := \qbhe/L\ [\unit{\watt\per\meter}]$, $\mathbf{x} = (\mathbf{x}_{xy}, z)$ with $\mathbf{x}_{xy}=(x,y)$, and $\mathbf p = (\mathbf{p}^1, \dots, \mathbf{p}^{\Nb})$ collects the planar positions $\mathbf{p}^k = (p_x^k, p_y^k)$ of the BHEs, all of length $L\ [\unit \meter]$.

The physical subsurface is the half-space $\mathbb{R}^3_+ := \{\mathbf{x} \in \mathbb{R}^3 : z \ge 0\}$, whose surface is held at $T_{\text{surface}} = T_0(0)$ \cite[Sec.~2.2]{CLAESSON1988509}, i.e.\ $u|_{z=0}=0$. We enforce this by the method of images: the second indicator function in \eqref{eq: source f} mirrors the BHE to $z<0$, so $f$ is odd in $z$. Since $v_z = 0$, $u$ is odd in $z$ as well. Hence the restriction of $u$ to $\mathbb{R}^3_+$ solves the half-space problem.
Since $f$ is a measure, \eqref{eq: advection diffusion} is understood in the following distributional sense, where $\mathscr S'(\mathbb{R}^3)$ denotes the space of tempered distributions.

     \begin{definition}[Solution] \label{def: solution}
        Let $T>0$. A function $u\in C([0,T];\mathscr S'(\mathbb{R}^3))$ \emph{solves} \eqref{eq: advection diffusion} on $(0,T)$ if $u(\cdot,0)=0$ in $\mathscr S'(\mathbb{R}^3)$ and, for every $\phi\in C_c^\infty(\mathbb{R}^3)$ and every $\zeta\in C_c^\infty(0,T)$,
        \begin{equation} \label{eq: weak formulation}
            \int_0^T\!\big\langle u(t),\phi\big\rangle\,\zeta'(t)\,dt
            \;+\;\int_0^T\!\big\langle u(t),\,\mathbf{v}\cdot\nabla\phi+\alpha\Delta\phi\big\rangle\,\zeta(t)\,dt
            \;=\; -\frac{\alpha}{\lambda}\int_0^T\!\big\langle f(\cdot,t),\phi\big\rangle\,\zeta(t)\,dt .
        \end{equation}
    \end{definition}

The solution is expressed through the kernel
\begin{equation} \label{def: mathcal G}
    \mathcal{G}(r,z, \tau, L) := \Lambda(r, \tau) \, \mathcal{Z}(z, \tau, L), \qquad \Lambda(r, \tau) := C_\lambda \mathcal{R}(r, \tau) \mathcal{W}_{\mathbf{v}}(\tau),
\end{equation}
with $C_\lambda = 1/(4\pi \lambda)$, the radial diffusion kernel $\mathcal{R}(r, \tau) = \tau^{-1} e^{-r^2/(4\alpha\tau)}$ with $r(\mathbf{x}_{xy}; \mathbf{p}^k) = \|\mathbf{x}_{xy} - \mathbf{p}^k\|$, the advective decay factor $\mathcal{W}_{\mathbf{v}}(\tau) = e^{-\|\mathbf{v}\|^2 \tau/(4\alpha)}$, and the vertical finite source kernel
    \begin{equation} \label{eq: kernel Z}
        \mathcal{Z}(z, \tau, L) = \frac{1}{2} \left[ \text{erf}\left(\frac{L-z}{2\sqrt{\alpha\tau}}\right) + 2\text{erf}\left(\frac{z}{2\sqrt{\alpha\tau}}\right) - \text{erf}\left(\frac{L+z}{2\sqrt{\alpha\tau}}\right) \right],
    \end{equation}
which accounts for the finite length and the mirror image. The directional influence of the flow relative to the source position enters through $\mathcal{S}_{\mathbf{v}}(\mathbf{x}_{xy}, \mathbf{p}^k) = \exp( \langle \mathbf{x}_{xy} - \mathbf{p}^k, \mathbf{v} \rangle/(2\alpha))$. Both $\mathcal S_\mathbf{v}$ and $\mathcal W_\mathbf{v}$ result from the shift of the heat kernel by $\mathbf{v}\tau$ (cf.\ Lemma~\ref{lem: mfls representation}).

Let $\mu_k:=\delta_{p^k_x}\otimes\delta_{p^k_y}\otimes\big(\mathbb{1}_{[0,L]}-\mathbb{1}_{[-L,0]}\big)$ be the line source of the $k$-th BHE together with its mirror image, so that $f(\cdot,t)=-\qbheL(t)\,\mu$ with the finite measure $\mu:=\sum_{k=1}^{\Nb}\mu_k$ of total variation $2L\Nb$. The fundamental solution of the operator $\partial_t+\mathbf{v}\cdot\nabla-\alpha\Delta$ is the heat kernel shifted by $\mathbf{v}\tau$ \cite[Sec.~2.3.1]{Evans2010PDE},
\begin{equation} \label{eq: Phi}
    \Phi(\mathbf{x},\tau):=(4\pi\alpha\tau)^{-3/2}\exp\Big(-\frac{\|\mathbf{x}-\mathbf{v}\tau\|^2}{4\alpha\tau}\Big), \qquad \mathbf{x}\in\mathbb{R}^3,\ \tau>0.
\end{equation}

\begin{theorem}[Existence, uniqueness and representation] \label{thm: existence}
    Let $\qbheL\in L^1_{\mathrm{loc}}([0,\infty))$ and $T>0$. Then \eqref{eq: advection diffusion} has exactly one solution $u$ in the sense of Definition~\ref{def: solution}, given by Duhamel's formula
    \begin{equation} \label{eq: duhamel}
        u(\mathbf{x},t):=\frac{\alpha}{\lambda}\int_0^t \bigl(\Phi(\cdot,t-s)*f(\cdot,s)\bigr)(\mathbf{x})\,ds=-\frac{\alpha}{\lambda}\int_0^t\qbheL(s)\,\bigl(\Phi(\cdot,t-s)*\mu\bigr)(\mathbf{x})\,ds, \qquad t\in[0,T].
    \end{equation}
    For every $t\in[0,T]$, $u(\cdot,t)$ is an integrable function with $\|u(\cdot,t)\|_{L^1}\le\frac{\alpha}{\lambda}\,2L\Nb\int_0^t|\qbheL|$, and at every point $\mathbf{x}=(\mathbf{x}_{xy},z)$ with $\mathbf{x}_{xy}\neq\mathbf{p}^k$ for all $k=1,\dots,\Nb$ it is given by the superposition of $\Nb$ Moving Finite Line Sources (MFLS) \cite{MFLS, TYEGINGRAS2014354}\footnote{For functions $f,g$ on $[0,\infty)$ we write $(f\ast g)(t):=\int_0^t f(t-\tau)g(\tau)\,d\tau$ for the (causal) convolution.}
    \begin{equation} \label{eq: mfls}
        u(\mathbf{x}, t; L, \mathbf{p}) = -\sum_{k=1}^{\Nb} \mathcal{S}_\mathbf{v}(\mathbf{x}_{xy}, \mathbf{p}^k) \bigl( \qbheL \ast \mathcal{G}_k \bigr)(t), \qquad \mathcal{G}_k(\tau) := \mathcal{G}\bigl(r(\mathbf{x}_{xy}; \mathbf{p}^k), z, \tau, L\bigr).
    \end{equation}
\end{theorem}

The proof is given in \ref{appendix: proof thm1}.

%% file: sections/temporal_coupling.tex
The system outlet temperature is the mean of the outlet temperatures $T_{\text{out}}^k$ of the $\Nb$ BHEs from \eqref{eq:outlet_final},
\begin{equation} \label{eq:outlet_mean}
T_{\text{out,sys}}(t) = \frac{1}{\Nb} \sum_{k=1}^{\Nb} T_{\text{out}}^k(t).
\end{equation}
Since \eqref{eq: advection diffusion} is linear, $u=\sum_{\ell}u_{\ell}$, where the $\ell$-th summand $u_{\ell}$ of \eqref{eq: mfls} is the response to BHE $\ell$ alone. The wall temperature of BHE $k$ is therefore
\begin{align}
T_{s,k}(z, t) &= T_0(z) + \sum_{\ell=1}^{\Nb} \Delta T_{k, \ell}(z, t), \label{eq:Ts_sum}\\
\Delta T_{k, \ell}(z, t) &=  -\,\mathcal{S}_{\mathbf{v}}(\mathbf{p}^{k}, \mathbf{p}^{\ell}) \bigl( \qbheL \ast \mathcal{G}(r_{k, \ell}, z, \cdot, L) \bigr)(t), \label{eq:DeltaT_kk}
\end{align}
where $\Delta T_{k,\ell}=u_{\ell}(\mathbf p^k,\cdot,\cdot)$ and $r_{k,\ell}=\|\mathbf{p}^k - \mathbf{p}^{\ell}\|$ for $\ell\neq k$. For $\ell=k$, the line source is singular at $\mathbf p^k$, so $u_k$ is evaluated on the borehole wall, $r_{k,k}=r_b$. On the wall, we set the directional factor $\mathcal S_{\mathbf v}$ to $1$.

Substituting \eqref{eq:Ts_sum} into \eqref{eq:outlet_final} and averaging by \eqref{eq:outlet_mean} decomposes the system outlet temperature into three potentials,
\begin{equation} \label{eq:potential_decomposition}
T_{\text{out,sys}}(t) = \Psi_{\text{ground}} - \Psi_{\text{soil}}(t) - \Psi_{\text{power}}(t),
\end{equation}
with
\begin{align}
    \Psi_{\text{ground}} &= \psi(L)  \int_0^L T_0(\xi)\, w_L(\xi)\, d\xi, \label{Psi: ground}\\
    \Psi_{\text{soil}}(t) &= \frac{\psi(L)}{\Nb} \sum_{k=1}^{\Nb} \sum_{\ell=1}^{\Nb}  \mathcal{S}_{\mathbf{v}}(\mathbf{p}^{k}, \mathbf{p}^\ell) \int_0^L w_L(\xi) \bigl( \qbheL \ast \mathcal{G}(r_{k, \ell}, \xi, \cdot, L) \bigr)(t)\, d\xi, \label{Psi: soil}\\
    \Psi_{\text{power}}(t) &= (\psi(L) - 1) \frac{\qbhe(t)}{\Cf}. \label{Psi: power}
\end{align}
Since $\Delta T_{k,\ell}\le 0$ for $\qbheL\ge 0$, extraction lowers the wall temperature, so $\Psi_{\text{soil}}$ enters \eqref{eq:potential_decomposition} with a negative sign.

\subsubsection{Ground Potential} \label{subsec: ground}
The ground potential $\Psi_{\text{ground}}$ is the contribution of the undisturbed subsurface temperature $T_0$. Inserting $N_W$ from \eqref{sys: constants} and the linear profile \eqref{eq: T0 linear} into \eqref{Psi: ground}, substituting $s=L-\xi$ and integrating $e^{M_l s}$ and $s\,e^{M_l s}$ over $[0,L]$ gives
\begin{equation} \label{eq: Psi ground closed}
\Psi_{\text{ground}} = \chi(L) \sum_{l=1}^2 C_l\, U_l(L), \qquad U_l(L) := L \left( T_{\text{surface}}\, \phi_1(M_l L) + g_{\text{geo}} L\, \phi_2(M_l L) \right),
\end{equation}
with $\chi(L) := \psi(L)/D_W(L)$ and the normalized exponential functions\footnote{For $|z|<0.1$, \eqref{eq: phi functions} loses accuracy by cancellation. There we evaluate the Taylor series $\phi_k(z)=\sum_{j\ge0}z^j/(j+k)!$ truncated after $10$ terms, which is accurate to machine precision and gives $\phi_1(0)=1$, $\phi_2(0)=\tfrac12$. For $|z|\ge0.1$, we evaluate \eqref{eq: phi functions} with \texttt{expm1} \cite{Tang1992-bx}.}
\begin{equation} \label{eq: phi functions}
\phi_1(z) = \frac{e^z - 1}{z}, \qquad \phi_2(z) = \frac{e^z - 1 - z}{z^2}.
\end{equation}

\subsubsection{Soil and Power Potential}
The soil potential $\Psi_{\text{soil}}$ collects the thermal interaction within the field. Interchanging the integrals in \eqref{Psi: soil} reduces it to a single temporal convolution,
\begin{equation} \label{eq: Psi_soil}
\Psi_{\text{soil}}(t) = \psi(L) \bigl( \qbheL \ast \mathcal{H}_\Sigma \bigr)(t), \qquad \mathcal{H}_\Sigma(\tau;\boldsymbol{\theta}) := \frac{1}{\Nb} \sum_{k=1}^{\Nb} \sum_{\ell=1}^{\Nb}  \mathcal{S}_{\mathbf{v}}(\mathbf{p}^{k}, \mathbf{p}^\ell)\, \mathcal{H}(r_{k, \ell}, \tau, L),
\end{equation}
with $\mathcal{H}(r, \tau, L):=\Lambda(r, \tau)\, \mathcal{Y}(\tau;L)$ and the integrated vertical response
\begin{equation} \label{def: mathcalV}
\mathcal{Y}(\tau;L) := \int_0^L \mathcal{Z}(z, \tau, L)\, w_L(z) \, dz,
\end{equation}
whose closed form is derived in \ref{appendix: closed V}.
The power potential $\Psi_{\text{power}}$ is the quasi-steady temperature drop across the borehole resistance; it does not depend on the positions of the boreholes.

%% file: sections/spectral_decomposition.tex
The convolutions in \eqref{eq: Psi_soil} are expensive, because BHE fields are simulated over several decades (Section~\ref{sec: Constraints}) and their number grows quadratically with $\Nb$. To reduce their cost, we lump the load and quantify the resulting error.

\begin{figure}[htbp]
    \centering
    \includegraphics[width=0.65\linewidth]{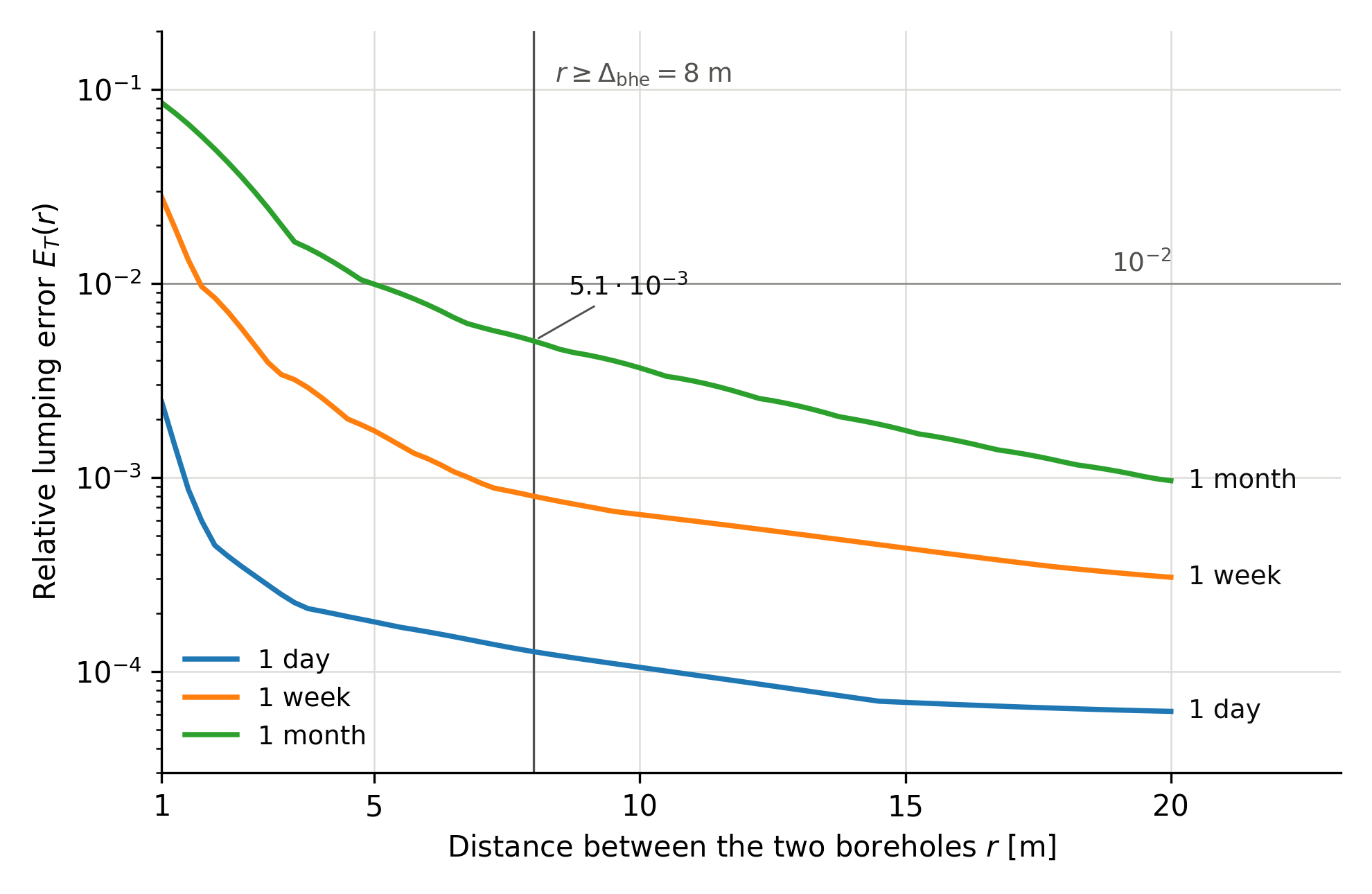}
    \caption{Lumping error $E_{\Delta t_B}(r)$ for lumping intervals $\Delta t_B$ of one day, one week and one month, as a function of the distance $r$ between two boreholes. Vertical line: minimum BHE distance $\Deltabhe=\SI{8}{\metre}$.}
    \label{fig: lumping error}
\end{figure}

Let $\Delta t_B$ be a lumping interval (e.g.\ a day, week or month), $t_m:=m\,\Delta t_B$, and let $\bar q_j$ be the mean of the load $q$ over the block $(t_{j-1},t_j]$. The convolution $\mathcal{C}(r,t) := (q \ast \mathcal{H}(r,\cdot))(t)$ is approximated at the block ends by
\begin{equation} \label{eq:lumping}
    \mathcal{C}_{\Delta t_B}(r,t_m) = \sum_{j=1}^{m}\bar q_j\int_{(m-j)\Delta t_B}^{(m-j+1)\Delta t_B}\mathcal{H}(r,\tau)\,d\tau ,
\end{equation}
and between the block ends $\mathcal{C}_{\Delta t_B}(r,\cdot)$ is interpolated linearly in $t$. For $\Delta t_B=1\,\mathrm{h}$, $\mathcal{C}_{1\,\mathrm{h}}=\mathcal{C}$ at the hourly time steps. We are interested in radii $r\ge\Deltabhe$, the minimum distance between two BHEs (Table~\ref{tab:parameters}). We measure the lumping error over the entire load history by $E_{\Delta t_B}(r) = \max_t | \mathcal{C}_{\Delta t_B}(r,t) - \mathcal{C}_{1\,\mathrm{h}}(r,t) | / \max_t | \mathcal{C}_{1\,\mathrm{h}}(r,t) |$.
Figure~\ref{fig: lumping error} shows $E_{\Delta t_B}(r)$ for the scaled case-study load (cf.\ Test 4 of \citet{ahmadfard2019review}).

For $\Delta t_B=1$ month the error is $8.5\,\%$ at $r=\SI{1}{\metre}$, exceeds $1\,\%$ for $r\lesssim\SI{5}{\metre}$, and is $0.51\,\%$ at $r=\Deltabhe$ and $0.17\,\%$ at $r=\SI{15}{\metre}$. For all $\SI{8}{\metre}\le r\le\SI{30}{\metre}$, it is at most $0.08\,\%$ for one week and $0.013\,\%$ for one day.
Since the error is largest at small radii, we split the soil potential into a self-interaction and a mutual interaction term, $\Psi_{\text{soil}}=\Psi_{\text{self}}+\Psi_{\text{inter}}$. The $\Nb$ diagonal terms $k=\ell$ of \eqref{eq: Psi_soil} are identical ($r_{k,k}=r_b$, $\mathcal{S}_\mathbf{v}=1$), so their sum cancels the factor $1/\Nb$ and gives $\Psi_{\text{self}}(t) = \psi(L) ( \qbheL \ast \mathcal{H}(r_b, \cdot, L) )(t)$. The remaining pairs give
\begin{equation} \label{eq: Psi_inter}
    \Psi_{\text{inter}}(t) = \frac{\psi(L)}{\Nb} \sum_{1 \leq \ell < k \leq \Nb} \hat{\mathcal{S}}_{\mathbf{v}}(\mathbf{p}^{\ell}, \mathbf{p}^k) \bigl( \qbheL \ast \mathcal{H}(r_{k, \ell}, \cdot, L) \bigr)(t),
\end{equation}
with the symmetrized shift kernel $\hat{\mathcal{S}}_\mathbf{v}(\mathbf{p}^{\ell}, \mathbf{p}^k) := \mathcal{S}_{\mathbf{v}}(\mathbf{p}^{\ell}, \mathbf{p}^k) + \mathcal{S}_{\mathbf{v}}(\mathbf{p}^{k}, \mathbf{p}^{\ell}) = 2 \cosh (\langle\mathbf{p}^k - \mathbf{p}^{\ell}, \mathbf{v} \rangle/(2 \alpha))$.
In $\Psi_{\text{inter}}$ all radii satisfy $r_{k, \ell}\ge\Deltabhe$.

%% file: sections/fft.tex
This subsection evaluates $u$ in \eqref{eq: mfls} and $\Psi_{\text{soil}}$ in \eqref{eq: Psi_soil} with the Discrete Fourier Transform (DFT) and derives their gradients, which the optimizer needs in every iteration.

\subsubsection{Evaluation of the Soil and Fluid Temperature}
We evaluate all convolutions on the hourly grid $t_n=n\Delta t$, $n=0,\dots,N$, and write $\mathcal{K}$ for either of the kernels $\mathcal{G}(r,z,\cdot,L)$ and $\mathcal{H}(r,\cdot,L)$. The heat load is piecewise constant, $\qbheL(s)=q_n$ for $s\in((n-1)\Delta t,\,n\Delta t]$, so splitting the integral at the grid points gives
\begin{equation} \label{eq:step convolution}
\int_0^{t_n}\qbheL(t_n-\tau)\,\mathcal{K}(\tau)\,d\tau=\sum_{j=1}^{n}q_{n-j+1}\,\kappa_j,\qquad \kappa_j:=\int_{(j-1)\Delta t}^{j\Delta t}\mathcal{K}(\tau)\,d\tau .
\end{equation}
For $r>0$, $\mathcal{K}$ and all its derivatives vanish as $\tau\to0^+$ because of the factor $e^{-r^2/(4\alpha\tau)}$ in $\mathcal{R}$. We compute the weights $\kappa_j$, i.e.\ the increments of $t\mapsto\int_0^t\mathcal{K}(\tau)\,d\tau$, with the cumulative trapezoidal rule. Its nodes are the logarithmically graded points $\tau_i:=\tau_{\mathrm{lo}}e^{ih}$, $i=0,\dots,n_\tau-1$, with $h:=\ln(t_N/\tau_{\mathrm{lo}})/(n_\tau-1)$, refined by the hourly nodes $t_j$.

The right-hand side of \eqref{eq:step convolution} is a discrete convolution. With the sequences $\tilde q=(q_1,\dots,q_N)$ and $\tilde\kappa=(\kappa_1,\dots,\kappa_N)$ padded with zeros to a common length $M\ge2N-1$, all $N$ time levels follow from
\begin{equation} \label{eq:convolution}
\bigl(\qbheL\ast\mathcal{K}\bigr)_n:=\sum_{j=1}^{n}q_{n-j+1}\kappa_j=\operatorname{DFT}^{-1}\left\{ \operatorname{DFT}(\tilde q)\cdot \operatorname{DFT}(\tilde\kappa) \right\}_{n},\qquad n=1,\dots,N,
\end{equation}
which costs $O(N\log N)$ with the Fast Fourier Transform (FFT) instead of $O(N^2)$ \cite[Sec.~2.3]{MARCOTTE2008651}.

The temperature change $u(\mathbf{x},t_n)$ follows from \eqref{eq: mfls} with the discrete convolutions $(\qbheL\ast\mathcal{G}_k)_n$. The terms $\Psi_{\text{self}}$ and $\Psi_{\text{inter}}$ (Section~\ref{sec: lumping split}) are evaluated in the same way with the kernels $\mathcal{H}(r_b,\cdot,L)$ and $\mathcal{H}_{k,\ell}(\tau):=\mathcal{H}(r_{k,\ell},\tau,L)$. Since the cost of $\Psi_{\text{inter}}$ grows with the number of pairs, its hourly loads and weights are replaced by the block means and block integrals of \eqref{eq:lumping}.

\subsubsection{Analytical Gradients} \label{subsec: jacobians}
The design variables are $\boldsymbol{\theta}=(L,\mathbf{p}^1,\dots,\mathbf{p}^{\Nb})$. The loads $\qsys$ and $\qbhe=\qsys/\Nb$ do not depend on $\boldsymbol{\theta}$, whereas $\qbheL=\qbhe/L$ depends on $L$. All gradients follow from differentiation under the integral sign \cite[Thm.~2.27]{Folland1999RA}. Let $q\in L^1(0,T)$, $\Theta\subseteq\mathbb{R}^d$ open, and let $K:(0,T]\times\Theta\to\mathbb{R}$ have first derivatives $D_{\boldsymbol{\theta}}K$ that are continuous on $(0,T]\times\Theta$ and satisfy
\begin{equation} \label{eq:domination}
    \sup_{\tau\in(0,T],\,\boldsymbol{\theta}\in\Theta_0}\bigl(|K(\tau;\boldsymbol{\theta})|+\|D_{\boldsymbol{\theta}}K(\tau;\boldsymbol{\theta})\|\bigr)<\infty \qquad\text{for every compact }\Theta_0\subset\Theta .
\end{equation}
Then, for every $t\in(0,T]$, $u(t;\boldsymbol{\theta}):=\bigl(q\ast K(\cdot;\boldsymbol{\theta})\bigr)(t)$ is continuously differentiable on $\Theta$ \cite[Thm.~9.21]{Rudin1976}, with
\begin{equation} \label{eq: derivative convolution}
D_{\boldsymbol{\theta}}u(t;\boldsymbol{\theta})=\bigl(q\ast D_{\boldsymbol{\theta}}K(\cdot;\boldsymbol{\theta})\bigr)(t).
\end{equation}
We apply this to $u$, $\Psi_{\text{soil}}$ and the outlet temperature $\Toutsys$; the proofs are given in \ref{Appendix_Jacobians}.

The following propositions use $\mathcal{G}_{\text{mid}}(r,\tau,L):=\Lambda(r,\tau)\Zmid(\tau,L)$ with $\Zmid(\tau,L):=\mathcal{Z}(L/2,\tau,L)$, i.e.\ the kernel for the soil temperature at mid-depth (cf.\ Section~\ref{sec: Constraints}), and $\mathbf{v}_{xy}=(v_x,v_y)$.

\begin{proposition}[Gradient of the soil temperature] \label{prop:soil jacobian}
Let $\mathbf{x}_{xy}\in\mathbb{R}^2$ be a sensor position, $\qbhe\in L^1(0,T)$ and $\Theta_x:=\{(L,\mathbf{p}):L>0,\ \mathbf{p}^k\neq\mathbf{x}_{xy}\ \text{for all }k\}$. For $\boldsymbol{\theta}\in\Theta_x$ let $u(t;\boldsymbol{\theta})$ be given by \eqref{eq: mfls} at $\mathbf{x}=(\mathbf{x}_{xy},L/2)$ with $\qbheL=\qbhe/L$, and let $r_k=\|\mathbf{x}_{xy}-\mathbf{p}^k\|$, $\mathcal{S}_k=\mathcal{S}_{\mathbf{v}}(\mathbf{x}_{xy},\mathbf{p}^k)$, $\mathcal{G}_k(\tau)=\mathcal{G}_{\text{mid}}(r_k,\tau,L)$. Then $u(t;\cdot)\in C^1(\Theta_x)$ for every $t\in(0,T]$, and
\begin{align}
\partial_L u(t;\boldsymbol{\theta}) &= -\frac{u(t;\boldsymbol{\theta})}{L}-\sum_{k=1}^{\Nb}\mathcal{S}_k\Bigl(\qbheL\ast\bigl[\Lambda(r_k,\cdot)\,\partial_L\Zmid(\cdot,L)\bigr]\Bigr)(t), \label{eq: dL u}\\
\nabla_{\mathbf{p}^k}u(t;\boldsymbol{\theta}) &= -\mathcal{S}_k\bigl(\qbheL\ast[\mathbf{b}_k\,\mathcal{G}_k]\bigr)(t), \qquad \mathbf{b}_k(\tau):=\frac{1}{2\alpha}\left[\frac{\mathbf{x}_{xy}-\mathbf{p}^k}{\tau}-\mathbf{v}_{xy}\right], \label{eq: dp u}
\end{align}
where $\partial_L\Zmid(\tau,L)=\dfrac{3}{4\sqrt{\pi\alpha\tau}}\left(e^{-\frac{L^2}{16\alpha\tau}}-e^{-\frac{9L^2}{16\alpha\tau}}\right)$.
\end{proposition}

\begin{proposition}[Gradient of $\Psi_{\text{soil}}$] \label{prop:psi jacobian}
Let $\qbhe\in L^1(0,T)$, $\Theta_\Psi:=\{(L,\mathbf{p}):L>0,\ \mathbf{p}^k\neq\mathbf{p}^\ell\ \text{for }k\neq\ell\}$, and let $\Psi_{\text{soil}}(t;\boldsymbol{\theta})$ be given by \eqref{eq: Psi_soil} with $\qbheL=\qbhe/L$ and $r_{k,k}=r_b$. Then $\Psi_{\text{soil}}(t;\cdot)\in C^1(\Theta_\Psi)$ for every $t\in(0,T]$, and
\begin{align}
\partial_L\Psi_{\text{soil}}(t;\boldsymbol{\theta}) &= \left(\frac{\psi'(L)}{\psi(L)}-\frac1L\right)\Psi_{\text{soil}}(t;\boldsymbol{\theta})
+\frac{\psi(L)}{\Nb}\sum_{k,\ell=1}^{\Nb}\mathcal{S}_{\mathbf{v}}(\mathbf{p}^k,\mathbf{p}^\ell)\Bigl(\qbheL\ast\bigl[\Lambda(r_{k,\ell},\cdot)\,\partial_L\mathcal{Y}(\cdot;L)\bigr]\Bigr)(t), \label{eq: dL Psi}\\
\nabla_{\mathbf{p}^m}\Psi_{\text{soil}}(t;\boldsymbol{\theta}) &= \frac{\psi(L)}{\alpha\,\Nb}\sum_{\ell\neq m}\Bigl(\qbheL\ast\Bigl[\mathcal{H}_{m,\ell}\Bigl(\mathbf{v}_{xy}\sinh a_{m\ell}-\cosh(a_{m\ell})\,\frac{\mathbf{p}^m-\mathbf{p}^\ell}{\tau}\Bigr)\Bigr]\Bigr)(t), \label{eq: dp Psi}
\end{align}
with $a_{m\ell}:=\langle\mathbf{p}^m-\mathbf{p}^\ell,\mathbf{v}\rangle/(2\alpha)$, $\psi'(L)=-\gamma/(2\nu\sinh^2(\gamma L))$, and $\partial_L\mathcal{Y}$ is given in \eqref{eq: dV}.
\end{proposition}

\begin{proposition}[Gradient of the outlet temperature] \label{prop:fluid jacobian}
Under the assumptions of Proposition~\ref{prop:psi jacobian}, the outlet temperature $\Toutsys(t;\boldsymbol{\theta})$ in \eqref{eq:potential_decomposition} is continuously differentiable on $\Theta_\Psi$ for every $t\in(0,T]$, and
\begin{equation*}
\partial_L\Toutsys=\Psi_{\text{ground}}'(L)-\partial_L\Psi_{\text{soil}}-\psi'(L)\frac{\qbhe(t)}{\Cf},\qquad \nabla_{\mathbf{p}^m}\Toutsys=-\nabla_{\mathbf{p}^m}\Psi_{\text{soil}},
\end{equation*}
with $\partial_L\Psi_{\text{soil}}$, $\nabla_{\mathbf{p}^m}\Psi_{\text{soil}}$ from Proposition~\ref{prop:psi jacobian}, $\Psi_{\text{ground}}'(L)=\chi'(L)\sum_{l=1}^2C_lU_l(L)+\chi(L)\sum_{l=1}^2C_l\bigl(T_0(L)+M_lU_l(L)\bigr)$, $\chi=\psi/D_W$, $\chi'(L)=\bigl(\psi'(L)D_W(L)-\psi(L)D_W'(L)\bigr)/D_W(L)^2$ with $D_W'(L)=M_1D_1e^{M_1L}+M_2D_2e^{M_2L}$, and $U_l$ from \eqref{eq: Psi ground closed}.
\end{proposition}

%% file: sections/constraints.tex
The BHE field is subject to geometric, operational and environmental constraints \cite{VDI4640Blatt2}, \cite[Sec.~4]{TSAGARAKIS20202556}.

\subsubsection{Geometric and Spatial Constraints}
To comply with property law, each borehole position $\mathbf{p}^k = (p_x^k, p_y^k)$, $k = 1, \dots, \Nb$, must lie in the effective domain $\Omega_{\text{eff}} := \{ \mathbf{x} \in \Omega : \operatorname{dist}(\mathbf{x}, \partial\Omega) \geq \Deltaneigh \}$ of the property $\Omega$, with $\Deltaneigh>0$ \cite[Secs.~3.4 and~3.5]{TSAGARAKIS20202556}.
Technical guidelines \cite{VDI4640Blatt2} require $\|\mathbf{p}^k - \mathbf{p}^{\ell}\| \geq \Deltabhe$ for all $k \neq \ell$ to limit thermal interference and structural risks.
To replace this non-convex constraint by linear inequalities, we partition $\Omega$ into Voronoi cells $V_k$ and shrink each cell to $\hat V_k\subset V_k\cap\Omega_{\text{eff}}$ such that $\operatorname{dist}(\hat V_k, \hat V_{\ell}) \geq \Deltabhe$ for $k \neq \ell$ (Sections~\ref{section: spatial Distribution via CVT} and~\ref{subsec: optimal sharpness}). Then $\mathbf{p}^k \in \hat V_k$ for all $k$ guarantees both the separation and the distance to the property boundary. This requires a convex polygon $\Omega$, as in the case study, so that the cells $V_k$ and $\Omega_{\text{eff}}$ are convex polygons. For a non-convex property, the cells can be non-convex or even disconnected; this extension is left for future work.

The borehole length is bounded, $L \in \mathcal L :=[L_{\min}, L_{\max}]$. In Germany, boreholes deeper than $100\,\mathrm{m}$ fall under \S~127 of the Federal Mining Act \cite{BBergG127}.

\subsubsection{Operational Thermal Constraints} \label{subsec: fluid constraints}
To keep the coefficient of performance of the heat pump high and to prevent freezing of the ground and the grout, the fluid temperature must stay above a safety threshold \cite{VDI4640Blatt2}. For high-temperature heat pumps, higher return temperatures from the field also reduce the required temperature lift and thus the energy consumption \cite{BARCOBURGOS2022122378}.
The outlet temperature $\Toutsys$ from~\eqref{eq:outlet_mean} must satisfy
\begin{equation} \label{constr: fluid}
T_{\max}^{\text{fluid}}\geq \Toutsys(t)  \geq T_{\min}^{\text{fluid}} \quad \forall t \in [0, t_{\mathrm{end}}].
\end{equation}
$T_{\min}^{\text{fluid}}$ is typically set between $-3^\circ\text{C}$ and $0^\circ\text{C}$ to protect the borehole grout \cite{VDI4640Blatt2}, while $T_{\max}^{\text{fluid}}$ depends on the heat pump. Both limits refer to the outlet temperature $\Toutsys$, i.e.\ the fluid returning from the field to the heat pump. The inlet temperature $\Tin(t)=\Toutsys(t)-\qsys(t)/(\dot m_{\mathrm{tot}}c_p)$ is an affine function of $\Toutsys$ for the given load. A limit on $\Tin$, e.g.\ due to the pipe material, can therefore be imposed in the same way.

\begin{remark}
\label{rem:fixed_extrema}
Imposing \eqref{constr: fluid} at all $N=175\,200$ hourly time steps of $20$ years would add $N$ constraints and their gradients in every iteration. We therefore impose it only at the critical times $t^*_{\min} = \arg\min_{t\in[0,t_{\mathrm{end}}]} \Toutsys(t; \boldsymbol{\theta}_0)$ and $t^*_{\max} = \arg\max_{t\in[0,t_{\mathrm{end}}]} \Toutsys(t; \boldsymbol{\theta}_0)$ of the initial design $\boldsymbol{\theta}_0$, i.e.\ $\mathcal{T}^*_{\min}=\{t^*_{\min}\}$ and $\mathcal{T}^*_{\max}=\{t^*_{\max}\}$ in $(\mathrm{P})$. These times may shift during the optimization. The final design is therefore checked against the full hourly time series of $\Toutsys$, and the optimization is repeated with updated critical times if a limit is violated. Further critical times, e.g.\ one per heating and cooling season, can be added in the same way. The length-only baseline in Section~\ref{sec: case study} evaluates the extremes of the full time series.
\end{remark}

\subsubsection{Environmental and Regulatory Impact Constraints}
German environmental regulations require that the operation of a BHE field does not substantially alter the thermal state of neighboring properties \cite{hamburg_leitfaden,lbeg_leitfaden}. Unmanaged shallow geothermal systems interfere thermally and hydraulically and lose efficiency, as \citet{PEREGO2022102417} show for a densely developed urban area. At sensors $\mathbf{x}_i:=(\mathbf{x}_i^{xy},L/2)$, $i\in\mathcal{N}_{\mathrm{s}}$, at mid-depth, whose planar positions $\mathbf{x}_i^{xy}$ are typically located on the property boundary $\partial \Omega$, we therefore require
\begin{equation} \label{constr:soil}
|u(\mathbf{x}_i, t_{\mathrm{end}})| \leq \Delta T_{\max} \qquad \forall i\in\mathcal{N}_{\mathrm{s}},
\end{equation}
which the optimizer handles as the two smooth constraints $\pm u(\mathbf{x}_i, t_{\mathrm{end}}) \leq \Delta T_{\max}$. At mid-depth, the temperature change of a finite line source is close to its maximum over the depth.

The permissible change depends on the local administrative practice. The Hamburg guideline allows at most $\Delta T_{\max} = 1\,\mathrm{K}$ relative to the undisturbed temperature \cite[Sec.~5.6]{hamburg_leitfaden}, and the guideline of Lower Saxony asks for the isolines of $\pm1\,\mathrm{K}$ \cite[pp.~31--32]{lbeg_leitfaden}. Since $\Delta T_{\max}$ is an input of the optimization, any such limit can be prescribed; the case study uses $\Delta T_{\max} = 2\,\mathrm{K}$ (Table~\ref{tab:parameters}). The limit is checked by a simulation over a prescribed observation period. This period is $25$ years in Hamburg \cite[Sec.~8.6]{hamburg_leitfaden} and usually $50$ years in Lower Saxony, where the temporal extreme within this period is required \cite[p.~31]{lbeg_leitfaden}. We impose \eqref{constr:soil} only at $t_{\mathrm{end}}$, which neglects seasonal maxima that occur earlier. If the limit has to hold at all times, further time points can be added as in Remark~\ref{rem:fixed_extrema}.

%% file: sections/Problem_definitions.tex
The objective is to minimize the total drilling length $\Nb L$ under the constraints of the previous section. Since $\Nb$ is fixed, this is equivalent to minimizing $L$. The thermal demand $\qsys$ is prescribed, so every design meets it, and only the temperatures it causes are constrained. With the design variables $\boldsymbol{\theta} = (L, \mathbf{p}^1, \dots, \mathbf{p}^{\Nb})$, the problem reads
\begin{equation} \tag{$\mathrm{P}$} \label{eq: problem}
    \begin{aligned}
        \min_{\boldsymbol{\theta}} & \quad L \\
        \text{s.t.} & \quad \mathbf{p}^k \in \hat{V}_k, \quad \forall k=1, \dots, \Nb \quad && \text{(Spatial Separation)} \\
        & \quad  \Toutsys(t^* ; \boldsymbol{\theta}) \geq T_{\min}^{\text{fluid}}, \quad \forall t^* \in \mathcal{T}^*_{\min} \quad && \text{(Lower Limit Outlet Temperature)} \\
        & \quad  \Toutsys(t^* ; \boldsymbol{\theta}) \leq T_{\max}^{\text{fluid}}, \quad \forall t^* \in \mathcal{T}^*_{\max} \quad && \text{(Upper Limit Outlet Temperature)} \\
        & \quad |u(\mathbf{x}_i, t_{\mathrm{end}}; \boldsymbol{\theta})| \leq \Delta T_{\max}, \quad \forall i\in\mathcal{N}_{\mathrm{s}} \quad && \text{(Soil Temperature)} \\
        & \quad L \in \mathcal{L}. && \text{(Length Bounds)}
    \end{aligned}
\end{equation}
The cells $\hat{V}_k=\hat{V}_k(\eta_*)$ are computed before the optimization and kept fixed (Section~\ref{subsec: optimal sharpness}). We solve $(\mathrm{P})$ with the \texttt{trust-constr} method of SciPy \cite{2020SciPy-NMeth}, using the analytical gradients of Section~\ref{subsec: jacobians}.

%% file: sections/Spatial_strategy.tex
We construct the cells $\hat V_k=\hat V_k(\eta_*)$ of Section~\ref{sec: Constraints}, for the sharpness $\eta_*$ chosen in Section~\ref{subsec: optimal sharpness}, from a Centroidal Voronoi Tessellation (CVT) \cite{spatial_tessellations}, which spreads the BHEs over the property. The use of a CVT and of Lloyd's algorithm for the BHE placement was introduced in \cite{inproceed_1stpaper}, with a uniform density for the conduction-only case. Here, the density depends on the groundwater flow (Section~\ref{subsec: flow-density}), and the CVT serves as the initial layout and as the linear position constraints of the gradient-based optimization.

\subsubsection{CVT Optimization Problem}
For the convex polygon $\Omega \subset \mathbb{R}^2$ of the property and $\Nb$ BHEs, we seek generators $\mathbf{g} = \{\mathbf{g}_k\}_{k=1}^{\Nb} \subset \Omega$ and the partition of $\Omega$ into their Voronoi cells $V_k = \{ \mathbf{x} \in \Omega : \|\mathbf{x} - \mathbf{g}_k\| \le \|\mathbf{x} - \mathbf{g}_\ell\|\ \forall \ell \}$ that minimize
\begin{equation} \label{eq:objective}
\mathcal{J}(\mathbf{g}) := \sum_{k=1}^{\Nb} \int_{V_k} \rho(\mathbf{x}) \|\mathbf{g}_k - \mathbf{x}\|^2 \, d\mathbf{x}
\end{equation}
for a density $\rho > 0$. At a minimum, each generator is the weighted centroid $\mathbf{g}_k = \int_{V_k} \mathbf{x} \rho(\mathbf{x})\, d\mathbf{x} / \int_{V_k} \rho(\mathbf{x})\, d\mathbf{x}$ of its cell. Lloyd's algorithm \cite{Lloyd1982} approximates a CVT by alternately computing the Voronoi cells of the generators and moving the generators to these centroids.

\subsubsection{Flow-Dependent Density Functions} \label{subsec: flow-density}
\begin{figure}[ht]
    \centering
    \begin{subfigure}{0.45\textwidth}
        \centering
        \includegraphics[width=\linewidth]{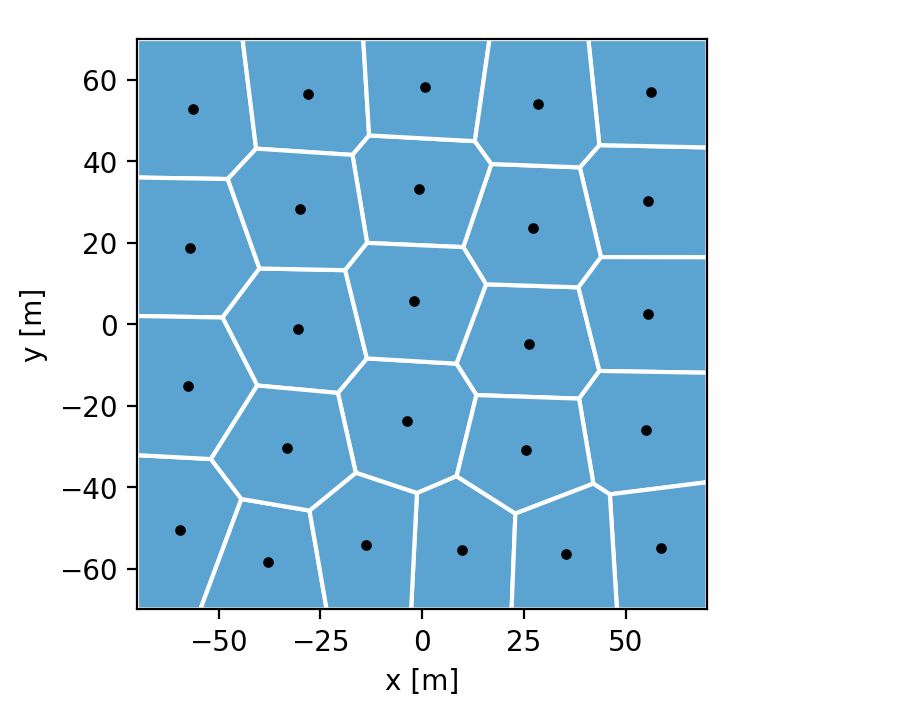}
        \caption{Uniform density, $\eta = 0$.}
    \end{subfigure}
    \hfill
    \begin{subfigure}{0.45\textwidth}
        \centering
        \includegraphics[width=\linewidth]{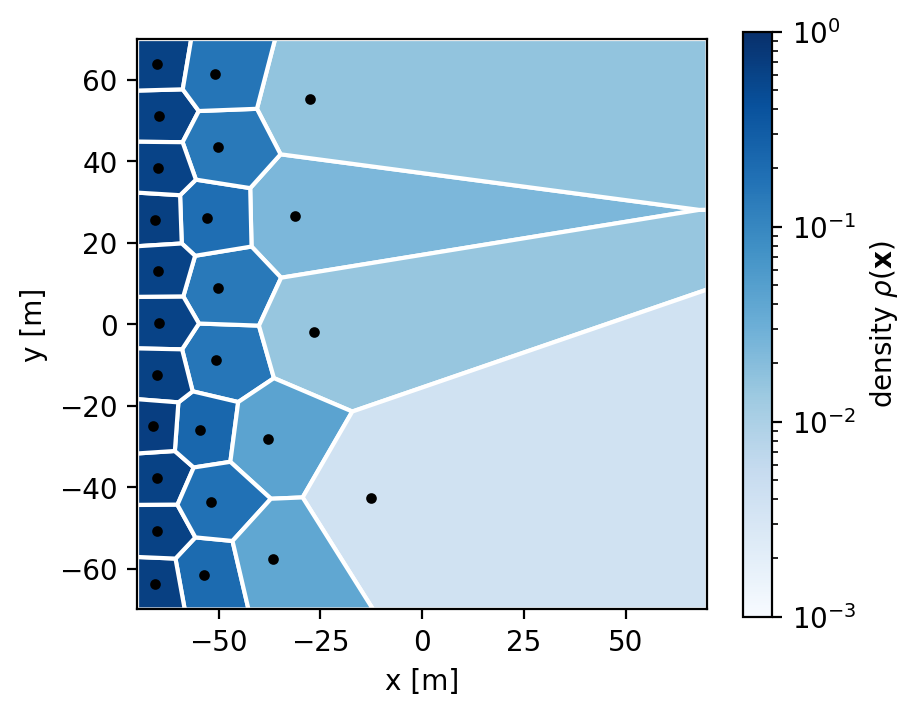}
        \caption{$\eta = \SI{0.0965}{\per\metre}$ (initial layout of the case study).}
    \end{subfigure}
    \caption{Centroidal Voronoi tessellations of the $\SI{140}{\metre}\times\SI{140}{\metre}$ property of Section~\ref{sec: case study} for $\Nb=25$ and groundwater flow in $x$-direction. The cells are colored by the density $\rho$ \eqref{eq: density rho} at their generator (black dots), normalized to $1$ at the upstream boundary (logarithmic scale).}
    \label{fig: voronoi non-uniform}
\end{figure}

With groundwater flow, upstream placement of the BHEs can reduce the thermal interference, because upstream boreholes are not affected by the plumes of other boreholes. We therefore use the density
\begin{equation} \label{eq: density rho}
\rho(\mathbf{x}) = \exp\left( -\eta\, \frac{\mathbf{x} \cdot \mathbf{v}}{\|\mathbf{v}\|} \right)
\end{equation}
in \eqref{eq:objective}, whose sharpness $\eta \ge 0$ controls the upstream bias. The density, and hence the generators $\mathbf{g}_k(\eta)$ and the cells $V_k(\eta)$, depend on $\eta$, which we indicate by the argument where it matters. A larger $\eta$ packs the cells more tightly upstream (Fig.~\ref{fig: voronoi non-uniform}). The resulting generators are the initial layout of the optimization.

\subsubsection{Analytical Centroid Calculation}
For the density \eqref{eq: density rho}, i.e.\ $\rho(\mathbf{x}) = e^{\mathbf{a} \cdot \mathbf{x}}$ with $\mathbf{a}=-\eta\,\mathbf{v}/\|\mathbf{v}\|$, the centroids in Lloyd's algorithm have a closed form: $\mathbf{g}_k(\eta)$ is the weighted centroid of $V_k(\eta)$ and follows from the following lemma with $V=V_k(\eta)$. We derive it with the divergence theorem in \ref{appendix: closed centroid}.

\begin{lemma} \label{lemma: centroid}
Let $V$ be a polygon in $\mathbb{R}^2$ with edges $S_i$ from vertex $\mathbf{P}_i$ to $\mathbf{P}_{i+1}$ (counter-clockwise). For the density $\rho(\mathbf{x}) = \exp(\mathbf{a} \cdot \mathbf{x})$ with $\mathbf{a} \neq 0$, the weighted centroid $\mathbf{g}$ is given by:
\begin{equation}
\mathbf{g} = \frac{\sum_{i} \left[ W_i\, \mathbf{m}_i - \Delta \mathbf{P}_i^\perp \mathcal{E}_{i,0} \right]}{\sum_{i} W_i \mathcal{E}_{i,0}}, \qquad \mathbf{m}_i := \mathbf{P}_i (\mathcal{E}_{i,0} - \mathcal{E}_{i,1}) + \mathbf{P}_{i+1} \mathcal{E}_{i,1},
\end{equation}
where $f_i = \mathbf{a} \cdot \mathbf{P}_i$, $\Delta f_i = f_{i+1} - f_i$, $\Delta \mathbf{P}_i^\perp = (y_{i+1}-y_i, -(x_{i+1}-x_i))^\top$ is the scaled outward normal, $W_i = \mathbf{a} \cdot \Delta \mathbf{P}_i^\perp$, $\mathcal{E}_{i,0} = e^{f_i}\, \phi_1(\Delta f_i)$ and $\mathcal{E}_{i,1} = e^{f_i} ( \phi_1(\Delta f_i) - \phi_2(\Delta f_i) )$, with $\phi_1, \phi_2$ from \eqref{eq: phi functions}, continuously extended by $\phi_1(0)=1$ and $\phi_2(0)=\tfrac12$.
\end{lemma}

\subsubsection{Choice of the Sharpness \texorpdfstring{$\eta$}{η}} \label{subsec: optimal sharpness}

For a given $\eta$, the weighted Lloyd iterations return the cells $V_k(\eta)$ and generators $\mathbf{g}_k(\eta)$. Each cell $V_k(\eta)=\{\mathbf{x}\mid \mathbf{n}_{kj}\cdot\mathbf{x}+b_{kj}\le0\ \forall j\}$, with the outward unit normals $\mathbf{n}_{kj}$ and offsets $b_{kj}$ of its edges $j$ (which depend on $\eta$ through the cell), is shrunk to $\hat{V}_k(\eta)=\{\mathbf{x}\mid \mathbf{n}_{kj}\cdot\mathbf{x}+b_{kj}+\varepsilon_{kj}\le0\ \forall j\}$, where $\varepsilon_{kj}=\varepsilon:=\Deltabhe/2$ for edges shared with a neighboring cell and $\varepsilon_{kj}=\Deltaneigh$ on $\partial\Omega$. A larger $\eta$ moves the generators $\mathbf{g}_k(\eta)$ closer to the upstream edges of their cells. We use the largest $\eta$ with $\mathbf{g}_k(\eta)\in\hat{V}_k(\eta)$ for all $k$, i.e.\ the strongest flow bias whose initial layout satisfies the separation constraint.

Algorithm~\ref{alg:maximal_beta} computes this $\eta$ for a fixed tessellation. Since the Lloyd iterations change the cells with $\eta$, its result is only an upper bound in the following procedure:
\begin{enumerate}
    \item Compute the uniform CVT ($\eta=0$), the starting point of all weighted Lloyd iterations.
    \item Apply Algorithm~\ref{alg:maximal_beta} to its cells $V_k(0)$, which gives the upper bound $\eta^{(0)}$.
    \item Run the weighted Lloyd iterations with $\eta=\eta^{(0)}$. If $\mathbf{g}_k(\eta^{(0)})\in\hat V_k(\eta^{(0)})$ for all $k$, accept $\eta^{(0)}$. Otherwise, bisect the interval $[0,\eta^{(0)}]$ for a fixed number of rounds; every trial value $\eta$ requires one Lloyd run, which gives $V_k(\eta)$, $\mathbf{g}_k(\eta)$ and $\hat V_k(\eta)$ for this $\eta$. Keep the largest $\eta$ with $\mathbf{g}_k(\eta)\in\hat V_k(\eta)$ for all $k$, or use $\eta=0$ if there is none. The accepted value is denoted $\eta_*$, and the cells $\hat V_k=\hat V_k(\eta_*)$ are kept fixed in the following.
\end{enumerate}

The functions \textsc{DistanceConstraint}$_k$ in Algorithm~\ref{alg:maximal_beta} need not be monotone, so the root returned by Brent's method \cite{Brent1973} is only a candidate. The feasibility of the final layout is checked in step~3.

\begin{algorithm}[H]
\caption{Largest admissible sharpness $\eta_{\max}^{\text{global}}$}
\label{alg:maximal_beta}
\begin{algorithmic}[1]
\Require Voronoi cells $\{V_k\}_{k=1}^{\Nb}$ (convex polygons), direction $\mathbf{v}$, margins $\varepsilon_{kj} > 0$, bracket upper bound $\eta_{\text{upper}}$
\State $\eta_{\max}^{\text{global}} \gets +\infty$

\For{each cell $V_k$ with vertices $\mathbf{P}_0, \dots, \mathbf{P}_{E_k-1}$ and edges $j=1,\dots,E_k$}
    \State Extract outward unit edge normals $\{\mathbf{n}_{kj}\}_j$, $\|\mathbf{n}_{kj}\|=1$, and offsets $\{b_{kj}\}_j$ such that $V_k = \{\mathbf{x} \mid \mathbf{n}_{kj} \cdot \mathbf{x} + b_{kj} \le 0\ \forall j\}$   
    \Function{DistanceConstraint$_k$}{$\eta$}
        \State Compute $\mathbf{g}_k(\eta)$ via Lemma \ref{lemma: centroid} ($V=V_k$, $\mathbf{a}=-\eta\mathbf{v}/\|\mathbf{v}\|$)
        \State \Return $\min_{j} \left( -b_{kj} - \mathbf{n}_{kj} \cdot \mathbf{g}_k(\eta) - \varepsilon_{kj} \right)$
    \EndFunction
    
    \If{$\text{DistanceConstraint}_k(0) \le 0$}
        \State $\eta_{\max, k} \gets 0$ \Comment{Uniform centroid already closer than $\varepsilon_{kj}$ to an edge}
    \ElsIf{$\text{DistanceConstraint}_k(\eta_{\text{upper}}) > 0$}
        \State $\eta_{\max, k} \gets \eta_{\text{upper}}$ \Comment{No sign change in the bracket}
    \Else
        \State Find a root $\eta_{\max, k}$ of $\text{DistanceConstraint}_k(\eta) = 0$ in $[0, \eta_{\text{upper}}]$ (Brent's method)
    \EndIf
    
    \State $\eta_{\max}^{\text{global}} \gets \min\left(\eta_{\max}^{\text{global}}, \eta_{\max, k}\right)$
\EndFor

\State \Return $\eta_{\max}^{\text{global}}$
\end{algorithmic}
\end{algorithm}

%% file: sections/forward_model_validation.tex
We compare the forward model with Earth Energy Designer (EED) \cite{etde_573620}, a commercial design tool, and the open-source GHEtool \cite{ghe_tool}. Both are based on the $g$-function method and are therefore not independent of each other, whereas OptiGeoS uses the MFLS with the temporal coupling above.

\begin{figure}[htbp]
    \centering
    \begin{subfigure}{0.9\linewidth}
        \centering
        \includegraphics[width=\linewidth]{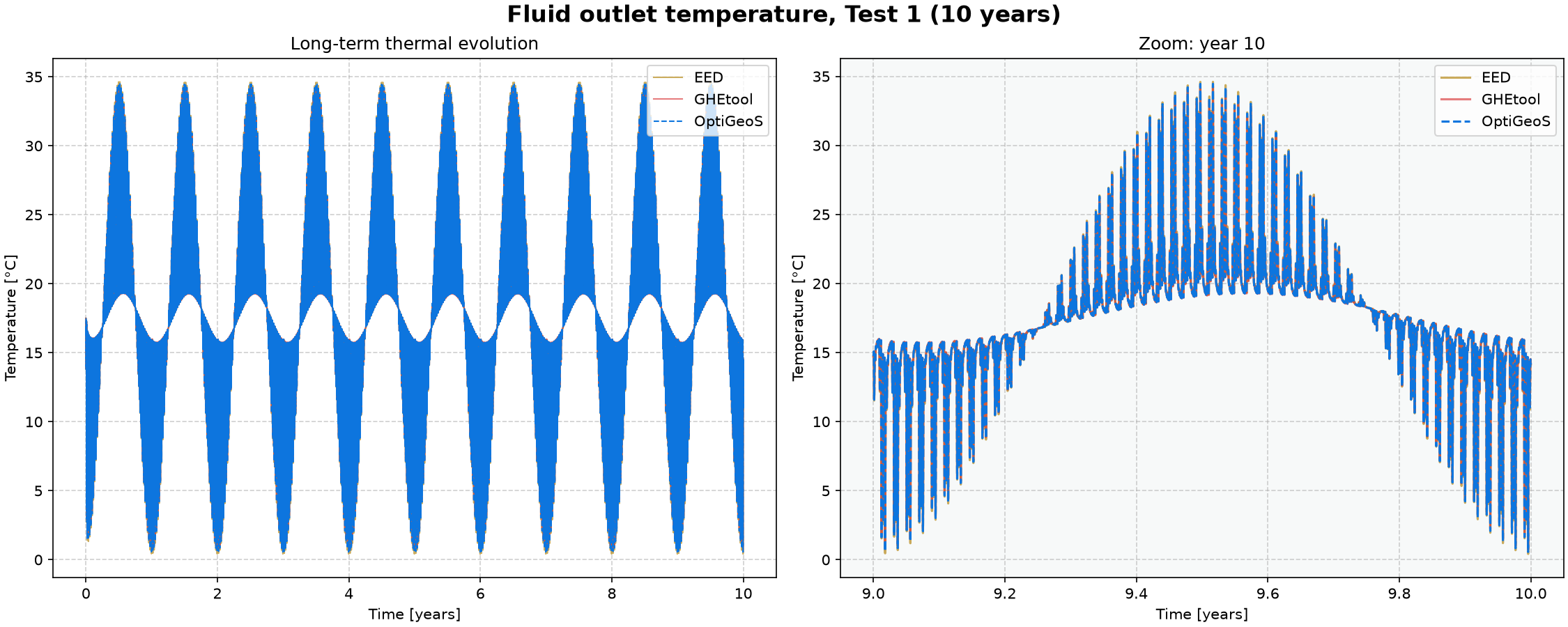}
        \caption{Test 1 (single borehole).}
    \end{subfigure}
    \begin{subfigure}{0.9\linewidth}
        \centering
        \includegraphics[width=\linewidth]{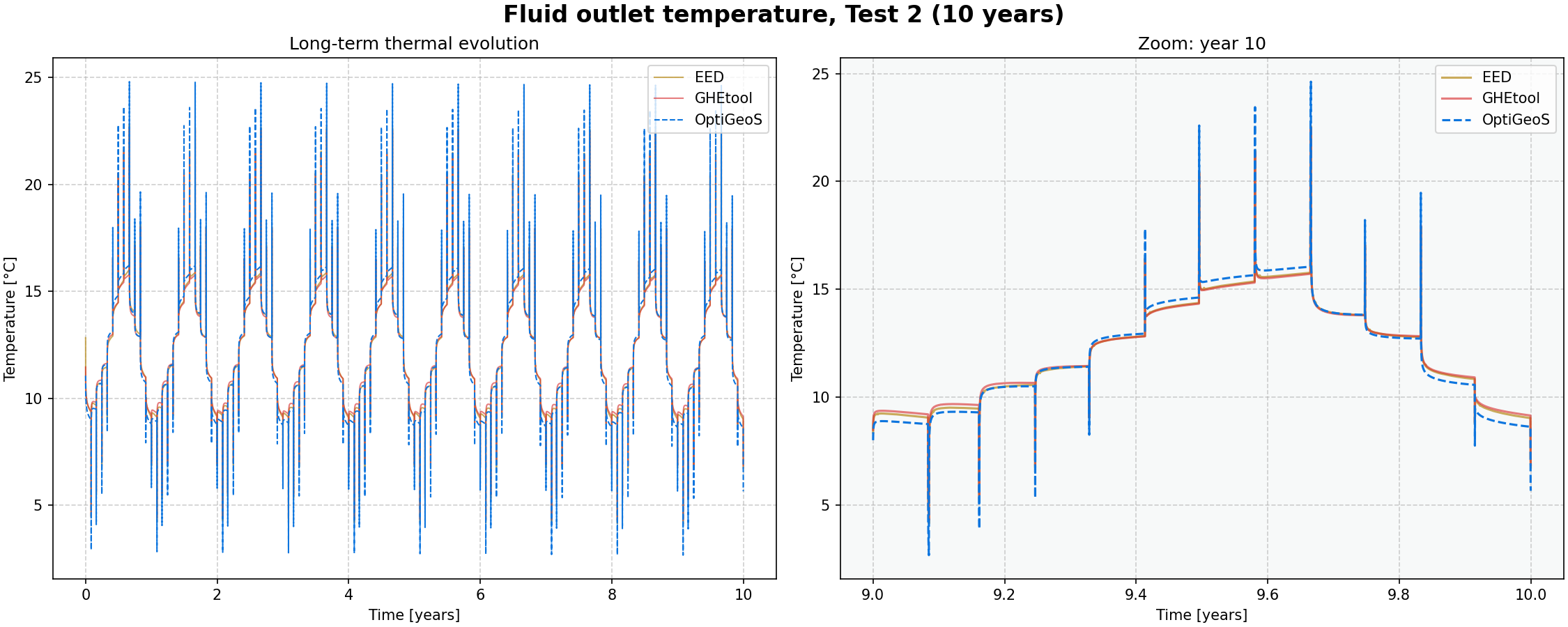}
        \caption{Test 2 ($120$ boreholes).}
    \end{subfigure}
    \caption{Fluid outlet temperature over $10$ years (left) and in the last year (right).}
    \label{fig:test results}
\end{figure}

We use Tests 1 and 2 of \citet[Secs.~4.2 and~4.3]{ahmadfard2019review}, which prescribe varying heat extraction rates. Test 1 is a single borehole of length $\SI{58.6}{\metre}$. Test 2 is a field of $120$ boreholes of length $\SI{90.1}{\metre}$ on a $12\times10$ grid with a spacing of $\SI{6}{\metre}$, so it also covers the thermal interaction between the boreholes. Figure~\ref{fig:test results} shows the fluid outlet temperatures of the three tools over $10$ years with hourly resolution.
EED and GHEtool model conduction only, so the comparison covers the model without groundwater flow ($\mathbf{v}=0$). The advective part, which drives the results of the case study, is verified separately below. Tests 1 and 2 are sizing tests \cite{ahmadfard2019review}, for which a comparison of the required borehole lengths would be more informative; we compare the temperature extremes for the prescribed fields.

\begin{table}[htbp]
    \centering
    \caption{Fluid outlet temperature extremes ($^\circ$C) and root mean square error (RMSE, K) of the hourly outlet temperature with respect to EED, 10-year simulation.}
    \label{tab:temperature_comparison}
    \begin{tabular}{llccc}
        \toprule
        \textbf{Case} & \textbf{Metric} & \textbf{EED} & \textbf{GHEtool} & \textbf{OptiGeoS} \\
        \midrule
        \multirow{3}{*}{Test 1} & Max Temp & 34.63 & 34.31 & 34.49 \\
                                & Min Temp & 0.42  & 0.74  & 0.56  \\
                                & RMSE to EED & --  & 0.09  & 0.04  \\
        \midrule
        \multirow{3}{*}{Test 2} & Max Temp & 22.83 & 22.60 & 24.81 \\
                                & Min Temp & 4.14  & 4.30  & 2.67  \\
                                & RMSE to EED & --  & 0.11  & 0.28  \\
        \bottomrule
    \end{tabular}
\end{table}

For Test 1, the temperature extremes of the three tools agree within $0.32\,\mathrm{K}$, the difference between EED and GHEtool (Table~\ref{tab:temperature_comparison}). OptiGeoS lies between them, $0.14\,\mathrm{K}$ from EED in both extremes, with an hourly RMSE of $0.04\,\mathrm{K}$.

For Test 2, the hourly time series agree in the mean (RMSE $0.28\,\mathrm{K}$ with respect to EED and $0.33\,\mathrm{K}$ with respect to GHEtool, mean differences $-0.01\,\mathrm{K}$ and $-0.04\,\mathrm{K}$; Table~\ref{tab:temperature_comparison}). OptiGeoS, however, predicts more pronounced extremes: its maximum is $2.0\,\mathrm{K}$ higher and its minimum $1.5\,\mathrm{K}$ lower than those of EED. The difference to EED is almost proportional to the instantaneous load (correlation coefficient $-0.99$). It is concentrated in the peak-load hours, with a mean absolute value of $1.6\,\mathrm{K}$ in the $1\,\%$ of the hours with the largest load and $0.2\,\mathrm{K}$ otherwise. This is consistent with an effective borehole thermal resistance about $0.04\,\mathrm{m\,K/W}$ above that of EED and GHEtool, compared with $R_b=\SI{0.11}{\metre\kelvin\per\watt}$ from the multipole method. Since the long-term mean agrees, the deviation points to the local borehole model (borehole resistances, quasi-steady description without grout capacity) rather than to the interaction between the boreholes. We did not isolate these effects. The higher peaks suggest that OptiGeoS would size longer boreholes for Test 2. Without measurements from a real BHE field, we cannot tell which tool is closer to the true response.

\begin{figure}[htbp]
    \centering
    \includegraphics[width=\linewidth]{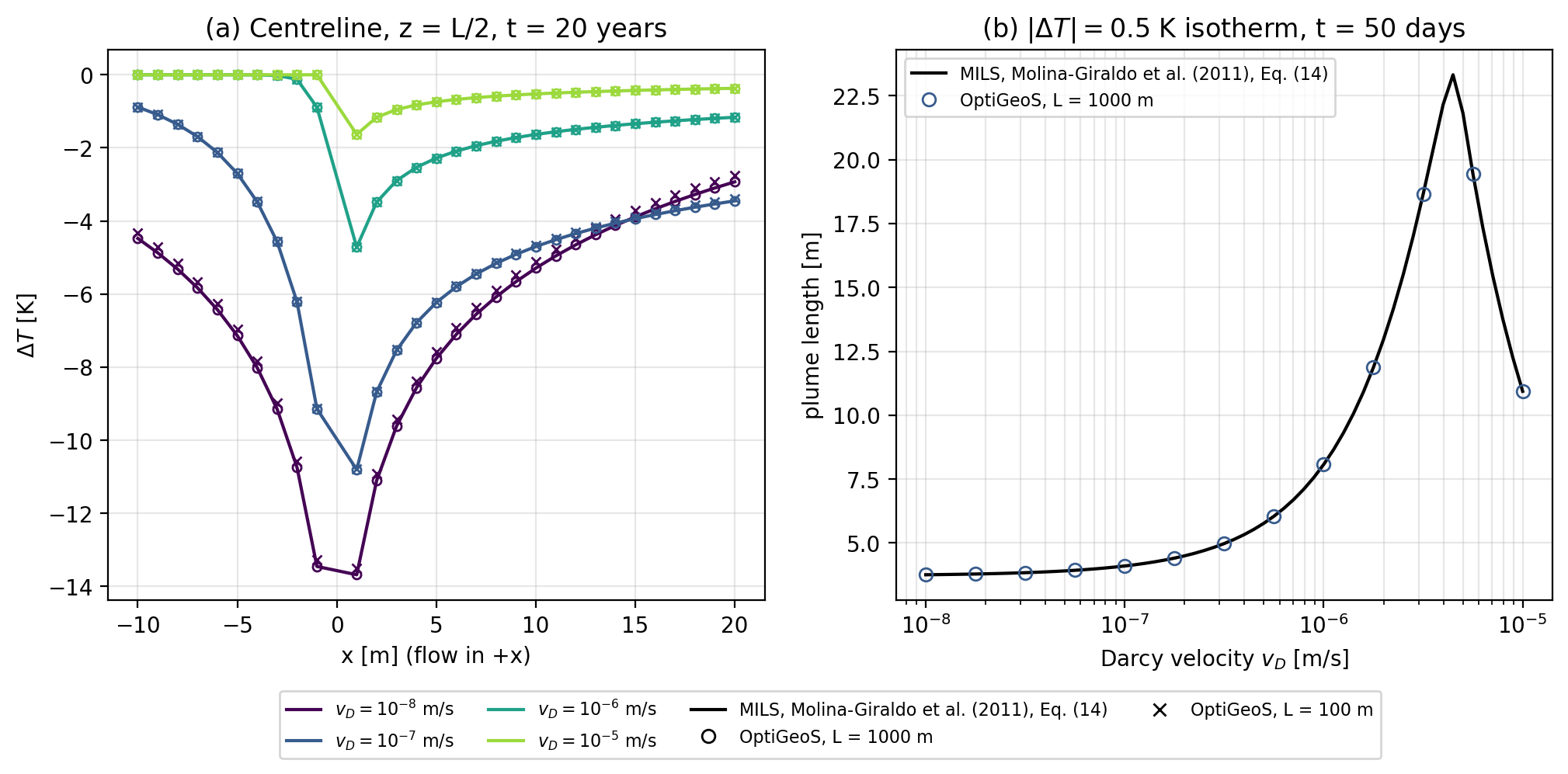}
    \caption{OptiGeoS compared with the MILS \cite[Eq.~(14)]{MolinaGiraldoetal2011} for a single borehole with $\SI{60}{\watt\per\metre}$ and Darcy velocity $\|\mathbf{v}_{\mathrm D}\|$. (a) Temperature change along the centerline at mid-depth after $20$ years. (b) Downstream distance of the $\SI{0.5}{\kelvin}$ isotherm after $50$ days, cf.\ \cite[Fig.~8]{MolinaGiraldoetal2011}.}
    \label{fig: mils}
\end{figure}

\paragraph{Advective soil response} To verify the advective part of the soil model, we compare it with the moving infinite line source (MILS) of \citet[Eq.~(14)]{MolinaGiraldoetal2011} without thermal dispersion. We use the setup of \citet[Sec.~4]{MolinaGiraldoetal2011}: a single borehole with a constant heat extraction rate of $\SI{60}{\watt\per\metre}$, $\lambda=\SI{2.5}{\watt\per\metre\per\kelvin}$, a volumetric heat capacity of $\SI{2.8e6}{\joule\per\cubic\metre\per\kelvin}$ and Darcy velocities $\|\mathbf{v}_{\mathrm D}\|$ between $10^{-8}$ and $\SI{e-5}{\metre\per\second}$. The effective thermal velocity is computed with the volumetric heat capacity of water $\SI{4.18e6}{\joule\per\cubic\metre\per\kelvin}$, which \citet{MolinaGiraldoetal2011} do not state. OptiGeoS evaluates the soil temperature at mid-depth with the same routine as the soil constraint. For a borehole of length $L=\SI{1000}{\metre}$, which approximates the infinite line source, the temperature change along the centerline ($-10\,\mathrm{m}\le x\le 20\,\mathrm{m}$) and perpendicular to the flow ($1\,\mathrm{m}\le y\le 10\,\mathrm{m}$) agrees with the MILS within $\SI{6e-5}{\kelvin}$ after $1$ and $20$ years for all velocities (Fig.~\ref{fig: mils}a shows $20$ years). The length of the $\SI{0.5}{\kelvin}$ plume after $50$ days agrees within $\SI{0.01}{\metre}$ (Fig.~\ref{fig: mils}b). For $L=\SI{100}{\metre}$, the finite length reduces the temperature change at mid-depth after $20$ years by at most $1.4\,\%$ ($\|\mathbf{v}_{\mathrm D}\|=\SI{e-8}{\metre\per\second}$), in line with the axial effects discussed in \cite{MFLS}. The monthly lumping of the interaction kernels (Section~\ref{sec: lumping split}) is accurate at these times. At $50$ days, it causes deviations of up to $\SI{0.15}{\kelvin}$, so Fig.~\ref{fig: mils}b is computed with daily blocks. This comparison is not a validation against numerical models of heterogeneous aquifers or against field measurements.

%% file: sections/case_study.tex
We compare a fixed-grid baseline with the simultaneous optimization of the positions and the length of a 1U BHE field under groundwater advection.

\subsection{Simulation Parameters and Constraints}
The geological and technical parameters follow Test 4 of \citet[Sec.~4.5 and Table~3]{ahmadfard2019review}. Table~\ref{tab:parameters} lists these parameters, the borehole geometry with the computed resistances, and the constraint bounds.
The field consists of $\Nb=25$ boreholes on a square property of $\SI{140}{\metre}\times\SI{140}{\metre}$. The thermal load is the load profile of Test 4, which was defined for a field of $120$ boreholes, scaled by the factor $\Nb/120$. The hourly load of one year ($8760$ values) is repeated over the simulated years. After the scaling, the peak cooling load of the field is $139.7\,\mathrm{kW}$ and the peak heating load $64.9\,\mathrm{kW}$. Over one year, $193.1\,\mathrm{MWh}$ of heat are injected into and $18.2\,\mathrm{MWh}$ extracted from the ground, so the load is cooling-dominated (cf.\ Fig.~\ref{fig: fluid_solution}). The simulated operation time is $t_{\mathrm{end}}=20$ years with an hourly time step, and $\Psi_{\text{inter}}$ is lumped with $\Delta t_B=730\,\text{h}$ (one month, Section~\ref{sec: lumping split}). The weights $\kappa_j$ in \eqref{eq:step convolution} are computed with $\tau_{\mathrm{lo}}=10^{-2}\,\mathrm{s}$ and $n_\tau=4000$, so that $h=6.2\cdot10^{-3}$.

\begin{table}[htbp]
\centering
\setstretch{1.0}
\small
\caption{Parameters of the case study: soil and fluid properties, borehole geometry with the resulting thermal resistances (1U), and constraint bounds.}
\label{tab:parameters}
\begin{tabular}{@{}llp{0.3cm}ll@{}}
\toprule
\multicolumn{2}{@{}l}{\textbf{Soil and fluid}} & & \multicolumn{2}{l@{}}{\textbf{Borehole (1U)}} \\
\midrule
Surface temperature $T_{\text{surface}}$ & \SI{15}{\celsius} & & Borehole radius $r_b$ & \SI{75}{\milli\metre} \\
Thermal conductivity $\lambda$ & \SI{1.9}{\watt\per\metre\per\kelvin} & & Shank spacing (centers) & \SI{83}{\milli\metre} \\
Thermal diffusivity $\alpha$ & \SI{9.26e-7}{\square\metre\per\second} & & Pipe outer diameter / wall & \SI{33.4}{} / \SI{3.7}{\milli\metre} \\
Geothermal gradient $g_{\text{geo}}$ & \SI{0}{\kelvin\per\metre} & & Conductivity grout / pipe & \SI{0.69}{} / \SI{0.4}{\watt\per\metre\per\kelvin} \\
Advection velocity $\mathbf{v}$ & \SI{1e-7}{\metre\per\second} in $x$ & & $R_b$ / $R_a$ (multipole) & \SI{0.2084}{} / \SI{0.8204}{\metre\kelvin\per\watt} \\
Fluid heat capacity $c_p$ & \SI{4019}{\joule\per\kilogram\per\kelvin} & & $R^\Delta_1 = R^\Delta_2$ / $R^\Delta_{12}$ & \SI{0.4167}{} / \SI{52.29}{\metre\kelvin\per\watt} \\
Total mass flow rate $\dot{m}_{\mathrm{tot}}$ & \SI{10.3401}{\kilogram\per\second} & & & \\
\cmidrule{4-5}
Fluid density & \SI{1026}{\kilogram\per\cubic\metre} & & \multicolumn{2}{l@{}}{\textbf{Constraints}} \\
Fluid thermal conductivity & \SI{0.468}{\watt\per\metre\per\kelvin} & & $\Deltabhe$ / $\Deltaneigh$ & \SI{8}{} / \SI{4}{\metre} \\
Fluid dynamic viscosity & \SI{3.37e-3}{\kilogram\per\metre\per\second} & & $T_{\min}^{\text{fluid}}$ / $T_{\max}^{\text{fluid}}$ & \SI{0}{} / \SI{38}{\celsius} \\
 & & & $\Delta T_{\max}$ & \SI{2}{\kelvin} \\
 & & & $L_{\min}$ / $L_{\max}$ & \SI{10}{} / \SI{150}{\metre} \\
\bottomrule
\end{tabular}
\end{table}

\subsection{Fixed-Grid Baseline}
As a baseline, we fix a regular $5\times5$ grid of BHEs with a spacing of $\SI{31}{\metre}$ and a distance of $\SI{8}{\metre}$ to the property boundary. Under groundwater flow, a designer would move the grid towards the inflow boundary. We therefore shift the grid by $\SI{4}{\metre}$ in the direction opposite to the flow, so that the upstream boreholes are $\SI{4}{\metre}$ from the upstream boundary. This is $(\mathrm{P})$ with fixed positions $\mathbf{p}^k$ and $L$ as the only variable. Let $\Gamma(L)$ denote the maximal violation of the fluid and soil temperature constraints for the length $L$. We compute a root of $\Gamma$ with Brent's method \cite{Brent1973}, which requires a bracket $\Gamma(L_{\min})>0>\Gamma(L_{\max})$. If $\Gamma(L_{\max})>0$, the baseline is infeasible within the drilling bounds; if $\Gamma(L_{\min})\le 0$, the lower bound $L_{\min}$ is already feasible. Since $\Gamma$ need not be monotone (e.g.\ the sensor depth $L/2$ changes with $L$), the root need not be the smallest feasible length.

\begin{figure}[htbp]
    \centering
    \includegraphics[width=0.8\linewidth]{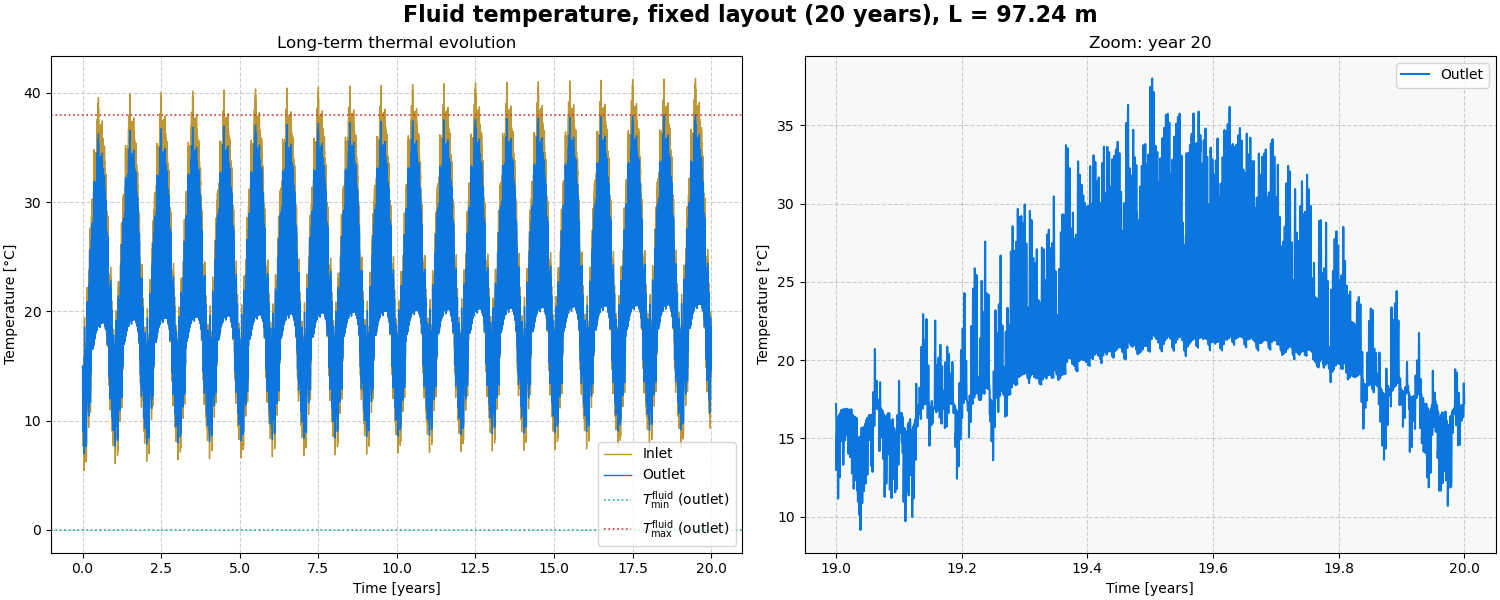}
    \caption{Inlet and outlet fluid temperature over 20 years for the fixed grid without groundwater flow ($L = \SI{97.24}{\metre}$). Dotted: limits $T_{\min}^{\text{fluid}}$ and $T_{\max}^{\text{fluid}}$ of the outlet temperature.}
    \label{fig: fluid_solution}
\end{figure}

Without groundwater flow ($\mathbf{v} = 0$), the baseline requires $L = \SI{97.24}{\metre}$, and the outlet temperature reaches its upper limit $T_{\max}^{\text{fluid}}$ in the last year of the 20-year period (Fig.~\ref{fig: fluid_solution}). At the peak cooling load, the inlet temperature exceeds the outlet temperature by the temperature spread of $\SI{3.4}{\kelvin}$ and reaches $\SI{41.3}{\celsius}$. This is admissible, since the limits refer to the outlet temperature (Section~\ref{subsec: fluid constraints}).

\begin{figure}[htbp]
    \centering
    \begin{subfigure}{0.45\textwidth}
        \centering
        \includegraphics[width=\linewidth]{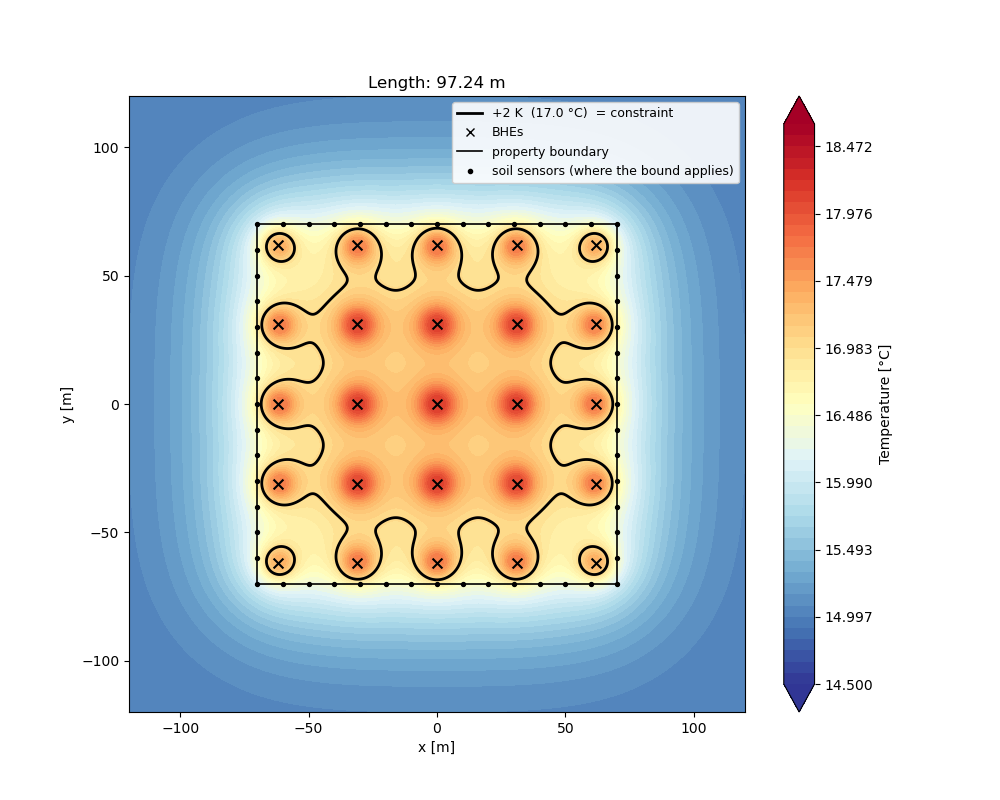}
        \caption{Fixed grid, $\mathbf{v}=0$.}
    \end{subfigure}
    \hfill
    \begin{subfigure}{0.45\textwidth}
        \centering
        \includegraphics[width=\linewidth]{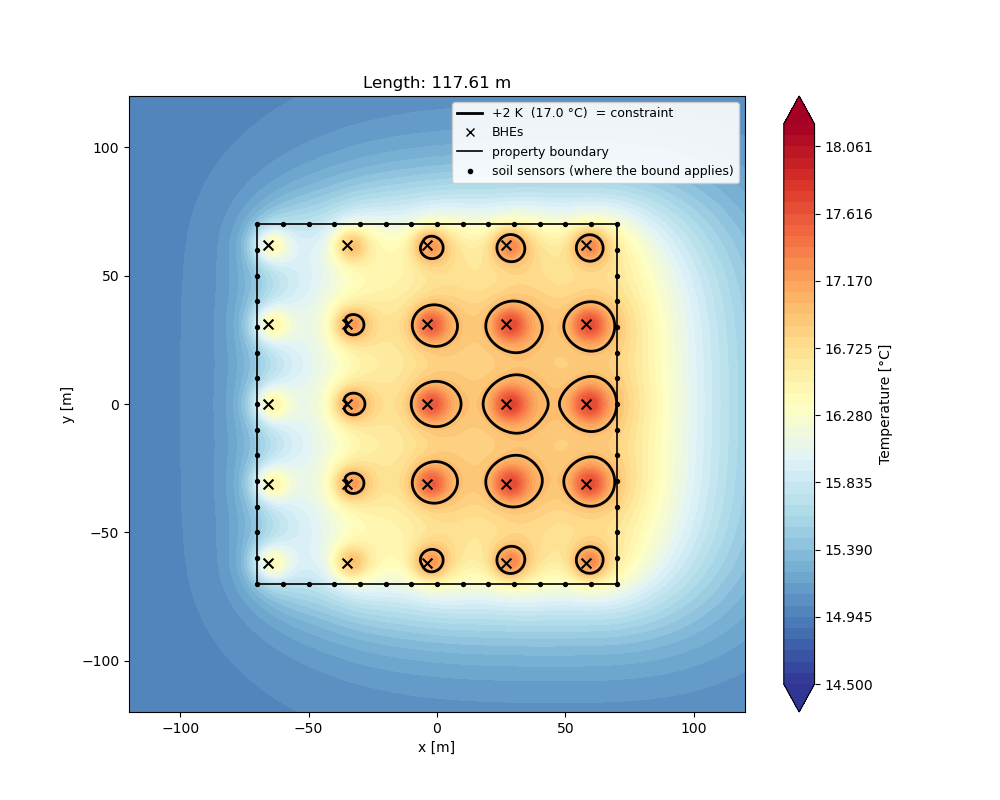}
        \caption{Fixed grid, $\|\mathbf{v}\|=\SI{1e-7}{\metre\per\second}$.}
    \end{subfigure}
    \caption{Soil temperature at mid-depth after $20$ years for the fixed $5\times5$ grid without groundwater flow and, shifted by $\SI{4}{\metre}$ towards the inflow boundary, with groundwater flow. $\times$: boreholes; black dots: sensors $\mathbf{x}_i$ of constraint \eqref{constr:soil}; black contour: $\Delta T_{\max}$.}
    \label{fig: MFLS case}
\end{figure}

With a groundwater flow of $\|\mathbf{v}\| = \SI{1e-7}{\metre\per\second}$ in $x$-direction, the thermal plume of the fixed grid is shifted downstream (Fig.~\ref{fig: MFLS case}). The soil constraint \eqref{constr:soil} becomes active at the downstream property boundary, although it is now $\SI{12}{\metre}$ from the outermost boreholes. The baseline then requires $L = \SI{117.61}{\metre}$. In this case study, the advective shift of the plume thus increases the required length by $\SI{20.37}{\metre}$ ($20.9\,\%$), which suggests that it should be considered in BHE field design.

\subsection{Optimized Layout}
For variable positions, we compute the CVT with the flow-dependent density of Section~\ref{subsec: flow-density}. Algorithm~\ref{alg:maximal_beta} and the subsequent bisection of Section~\ref{subsec: optimal sharpness} give $\eta = \SI{0.0965}{\per\metre}$.

\begin{figure}[htbp]
    \centering
    \begin{subfigure}{0.42\textwidth}
        \centering
        \includegraphics[width=\linewidth]{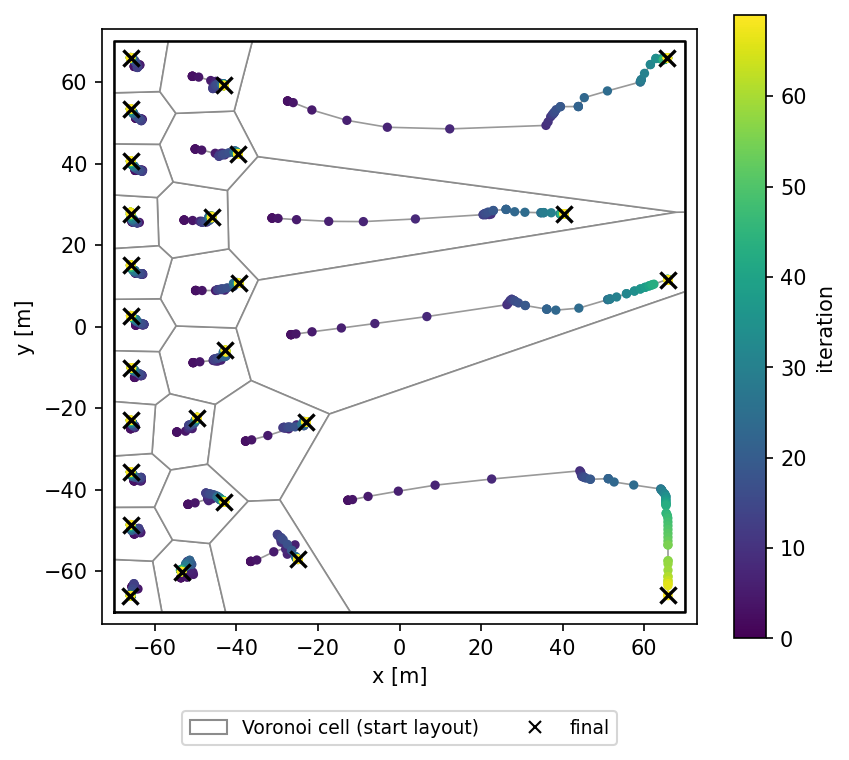}
        \caption{Borehole positions over the iterations.}
        \label{fig: optimization positions}
    \end{subfigure}
    \hfill
    \begin{subfigure}{0.42\textwidth}
        \centering
        \includegraphics[width=\linewidth]{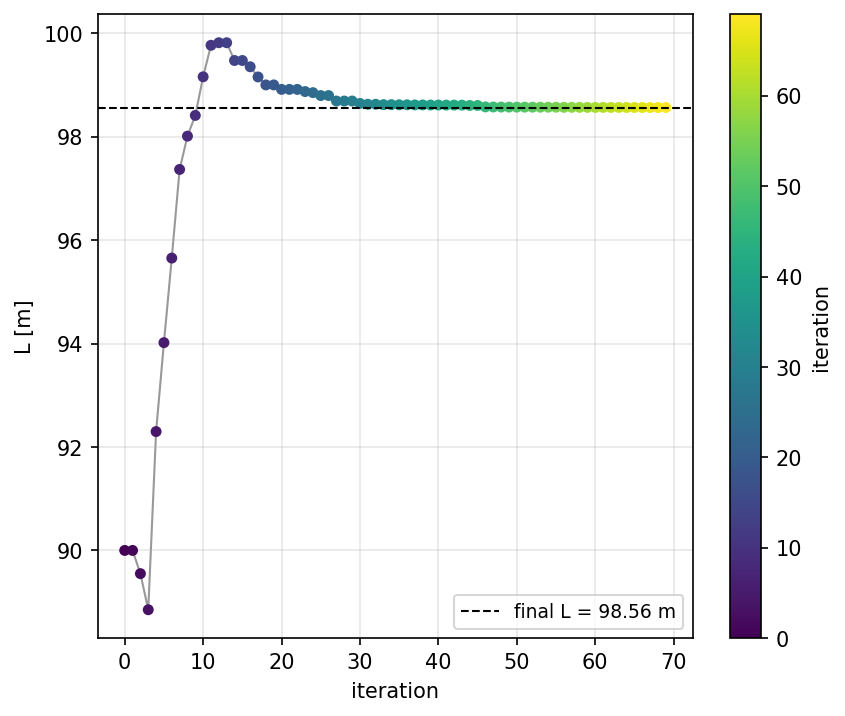}
        \caption{Borehole length $L$ over the iterations.}
        \label{fig: optimization length}
    \end{subfigure}
    \caption{Simultaneous optimization of the positions and the length under groundwater flow in $x$-direction, from the initial flow-biased CVT layout (dark) to the final iteration (bright). Gray lines in (a): Voronoi cells of the initial layout, each of which contains one borehole throughout the optimization.}
    \label{fig: optimization}
\end{figure}

\begin{figure}[htbp]
    \centering
    \includegraphics[width=0.6\linewidth]{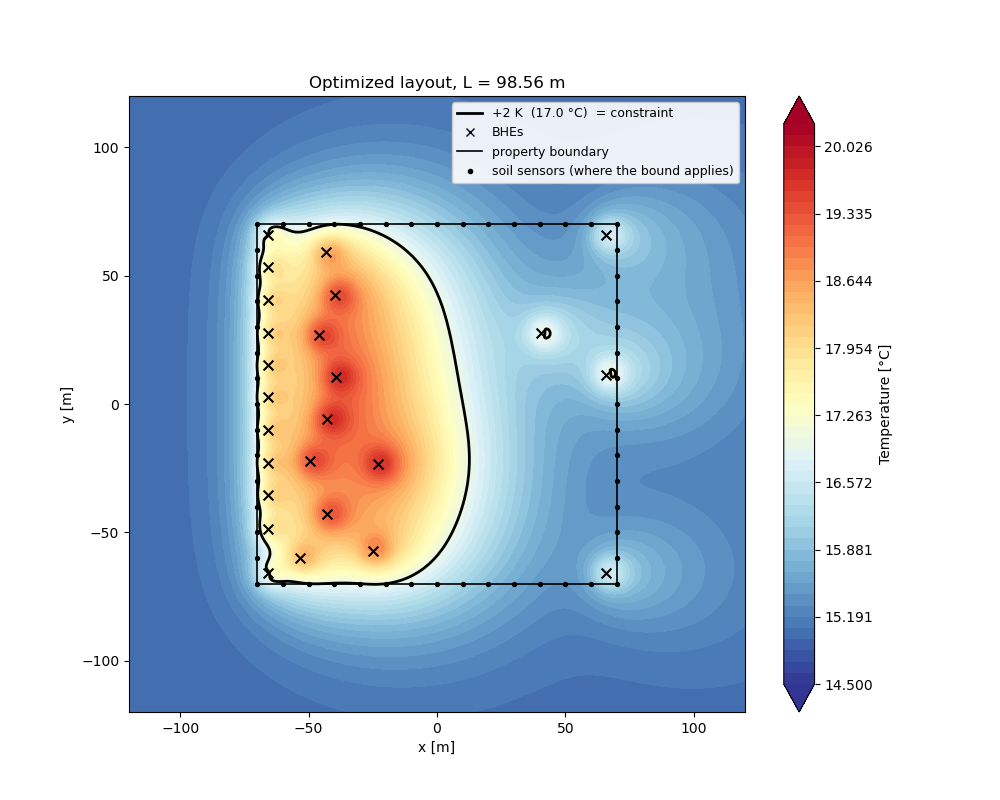}
    \caption{Soil temperature at mid-depth after $20$ years for the optimized layout ($L=\SI{98.56}{\metre}$) under groundwater flow; markers as in Fig.~\ref{fig: MFLS case}.}
    \label{fig: soil optimized}
\end{figure}

When the BHE positions adapt to the flow (Figs.~\ref{fig: optimization} and \ref{fig: soil optimized}), the required borehole length decreases from $\SI{117.61}{\metre}$ for the fixed grid to $L = \SI{98.56}{\metre}$, a reduction by $16.2\,\%$. At this optimum, two constraints are active: the soil constraint \eqref{constr:soil}, at sensors on the lateral property boundaries in the upstream half (Fig.~\ref{fig: soil optimized}), and the upper fluid temperature limit. Over the full hourly time series of $20$ years, the outlet temperature of the final layout stays between $7.1\,^\circ\mathrm{C}$ and $38.0\,^\circ\mathrm{C}$. The optimization starts from the infeasible length $L_0=\SI{90}{\metre}$ and increases $L$ to at most $\SI{99.82}{\metre}$ in iteration $12$ while the positions move. After iteration $30$, $L$ changes by less than $\SI{0.1}{\metre}$ (Fig.~\ref{fig: optimization length}).

To test the dependence on the initial layout, we repeated the optimization from five further random initializations of the Lloyd iterations. They give $\eta$ between $0.1002$ and $\SI{0.1093}{\per\metre}$ and final lengths between $98.42$ and $\SI{98.71}{\metre}$. These local optima differ by less than $0.3\,\%$ in length.

On a smaller property of $\SI{120}{\metre}\times\SI{120}{\metre}$ with the same load and a $5\times5$ grid of spacing $\SI{26}{\metre}$, the fixed grid, shifted by $\SI{4}{\metre}$ towards the inflow boundary as above, cannot satisfy the soil constraint within $L_{\max}$. At $L=\SI{150}{\metre}$, the temperature change at the downstream boundary still exceeds $\Delta T_{\max}$ by $\SI{0.12}{\kelvin}$. The optimized layout ($\eta=\SI{0.0907}{\per\metre}$) is feasible on this property with $L=\SI{100.81}{\metre}$ (Table~\ref{tab:ablation}).

\begin{table}[htbp]
\centering
\caption{Required borehole length $L$ for the different layouts of the case study. The change refers to the fixed grid on the same property with the same flow.}
\label{tab:ablation}
\small
\begin{tabular}{lllrr}
\toprule
Property & Flow & Layout & $L$ [m] & Change \\
\midrule
\multirow{4}{*}{$140\times140\,\mathrm{m}$} & \multirow{2}{*}{$\|\mathbf{v}\|=10^{-7}\,\unit{\metre\per\second}$} & Fixed $5\times5$ grid, shifted upstream & 117.61 & -- \\
 & & Positions and $L$ optimized & 98.56 & $-16.2\,\%$ \\
\cmidrule{2-5}
 & \multirow{2}{*}{$\mathbf{v}=0$} & Fixed $5\times5$ grid & 97.24 & -- \\
 & & Positions and $L$ optimized ($\eta=0$) & 96.81 & $-0.4\,\%$ \\
\midrule
\multirow{2}{*}{$120\times120\,\mathrm{m}$} & \multirow{2}{*}{$\|\mathbf{v}\|=10^{-7}\,\unit{\metre\per\second}$} & Fixed $5\times5$ grid & infeasible & -- \\
 & & Positions and $L$ optimized ($\eta=\SI{0.0907}{\per\metre}$) & 100.81 & -- \\
\bottomrule
\end{tabular}
\end{table}

\subsection{Computational Cost}
For $\Nb=25$ boreholes and $20$ years with hourly resolution, one simulation of the coupled model takes $0.3$--$0.4\,\mathrm{s}$ on a laptop CPU (Intel Core i7-8565U, $1.80\,\mathrm{GHz}$). One evaluation including all analytical gradients takes $5$--$6\,\mathrm{s}$. The optimization of $(\mathrm{P})$ needed $69$ iterations and $92$ evaluations of the coupled model and took about $8\,\mathrm{min}$. The cost of the pair interactions grows quadratically with the number of boreholes, and that of the FFT convolutions like $N\log N$ in the number $N$ of time steps.

\subsection{Discussion}
In the fixed grid, the thermal plume is advected downstream, heat accumulates at the downstream boreholes, and the soil constraint becomes active at the downstream property boundary (Fig.~\ref{fig: MFLS case}), consistent with \cite[Sec.~4.2]{TYEGINGRAS2014354}. With static positions, this can only be compensated by a longer borehole, which lowers the heat extraction rate $\qbheL$ per unit length.

Throughout the optimization, the Voronoi cells confine each borehole to its own region and thereby guarantee the minimal spacing. Within their cells, the optimized layout shows two mechanisms (Figs.~\ref{fig: optimization positions} and \ref{fig: soil optimized}). First, $11$ boreholes move to the upstream property boundary, at the distance $\Deltaneigh$. They exchange heat with groundwater that has not been disturbed by other boreholes, and their plumes are carried into the property. Second, four boreholes whose cells extend into the downstream half of the property move to the downstream ends of their cells, three of them to the downstream boundary. They are far from the other boreholes (Fig.~\ref{fig: optimization positions}), so the temperature change at the downstream boundary stays below $\Delta T_{\max}$.

With these mechanisms, the heat input is distributed over the property (Fig.~\ref{fig: soil optimized}). \citet{HechtMendez2013} found that the benefit of optimizing the loads decreases with the flow velocity, because the demand is fixed and the advective heat supply already levels the temperature anomalies. Here, the soil constraint at the property boundary makes the advected plume the limiting factor. The smaller property shows that the optimization can also make a field feasible for which the fixed $5\times5$ grid is infeasible within the drilling bounds.

%% file: sections/conclusion.tex
We presented a gradient-based framework that optimizes the borehole length and the borehole positions of a BHE field simultaneously. It couples an analytical forward model including groundwater advection with analytical gradients and a Voronoi-based partition of the property. The partition confines each borehole to its own cell and thereby guarantees the minimal spacing. In the case study with groundwater advection, the optimized layout reduced the required borehole length from $117.61\,\mathrm{m}$ to $98.56\,\mathrm{m}$, i.e.\ by $16.2\,\%$ compared with a fixed grid shifted towards the inflow boundary. This length is only $1.4\,\%$ above that of the fixed grid without groundwater flow. The improvement results from the gradient-based optimization of the positions within their Voronoi cells, which distributes the boreholes over the property. On a smaller property, the optimized layout was feasible although the fixed grid was not. These results suggest that, under groundwater flow, the layout of a BHE field should be adapted to the local flow conditions rather than fixed to a predefined grid.

The results are subject to several limitations. Every borehole carries the same extraction rate (Section~\ref{subsec: local bhe}), although under advection the upstream and downstream boreholes see different wall temperatures. The optimizer exploits these differences, so the reported savings should be read with this approximation in mind. The subsurface is homogeneous with a constant, horizontal effective velocity, so layered soils, heterogeneities or vertical flow may distort the plumes. The quasi-steady borehole model does not resolve sub-hourly transients. The fluid temperature constraints are evaluated at the two critical times of the initial design, which may shift during the optimization for strongly dynamic loads. All boreholes share one length $L$; individual lengths could give additional savings, in particular downstream. Finally, the forward model was compared with the conduction-only tools EED and GHEtool and with the analytical moving infinite line source. EED and GHEtool agree with it in the long-term response but give up to $2\,\mathrm{K}$ less pronounced peak temperatures for the field of $120$ boreholes. The forward model has not been validated against numerical models of heterogeneous aquifers or measurements from a real BHE field, and the benchmark comparison does not cover the required borehole lengths.

Future work includes a validation with measured data and the extension of the spatial strategy to non-convex properties, e.g.\ by a decomposition of the property into convex parts \cite{HERTEL198552}. The framework can also be extended to generic ground response functions \cite[Sec.~4]{TYEGINGRAS2014354}, \cite{CIMMINO2014641}. Another direction is the optimization of several interacting BHE fields to maximize the regional capacity while controlling their mutual thermal interference. An example is a brownfield scenario in which a new field is designed around an existing one.

%% file: sections/Appendix_proof_mfls.tex
This appendix proves Theorem~\ref{thm: existence}. Throughout, $T>0$ and $\Phi$ is given by \eqref{eq: Phi}. Lemma~\ref{lem: uniqueness} concerns the homogeneous problem and does not involve the source; the source \eqref{eq: source f} of the BHE field enters Lemmas~\ref{lem: duhamel solves} and~\ref{lem: mfls representation} and the proof of Theorem~\ref{thm: existence}. We write $f(\cdot,t)=-\qbheL(t)\,\mu$ with the finite measure $\mu$ of total variation $\|\mu\|=2L\Nb$ defined before \eqref{eq: Phi}. The mirror term in \eqref{eq: source f} makes $\mu$ odd in $z$. Lemmas~\ref{lem: duhamel solves} and~\ref{lem: uniqueness} and Theorem~\ref{thm: existence} therefore concern the full-space problem with the source and its mirror image. Since $v_z=0$, its solution is odd in $z$, and its restriction to $z\ge0$ is the solution of the half-space problem with $u|_{z=0}=0$.

\begin{lemma}[Duhamel's formula gives a solution] \label{lem: duhamel solves}
    Let $\qbheL\in L^1(0,T)$ and let $f$ be the source \eqref{eq: source f}. Then, for every $t\in[0,T]$, the function $u(\cdot,t)$ defined by \eqref{eq: duhamel} is integrable with $\|u(\cdot,t)\|_{L^1}\le\frac{\alpha}{\lambda}\,2L\Nb\int_0^t|\qbheL|$, and $u$ solves \eqref{eq: advection diffusion} in the sense of Definition~\ref{def: solution}.
\end{lemma}
\begin{proof}
    Since $f(\cdot,s)=-\qbheL(s)\,\mu$ with a fixed finite measure $\mu$ and $\qbheL\in L^1(0,T)$, all integrands below are products of the measurable function $\qbheL$ and a bounded Borel function of $(s,\mathbf{y})$ (e.g.\ $(s,\mathbf{y})\mapsto(P_{t-s}\phi)(\mathbf{y})$ on $\{s<t\}$), so that Tonelli's and Fubini's theorems apply. Since $\Phi(\cdot,\tau)$ is a probability density, $|\Phi(\cdot,\tau)*f(\cdot,s)|\le|\qbheL(s)|\,\Phi(\cdot,\tau)*|\mu|$ and $\int_{\mathbb{R}^3}\Phi(\cdot,\tau)*|\mu|\,d\mathbf{x}=\|\mu\|=2L\Nb$. By Tonelli's theorem, \eqref{eq: duhamel} therefore defines $u(\cdot,t)\in L^1(\mathbb{R}^3)\subset\mathscr S'(\mathbb{R}^3)$ with the stated bound, and by Fubini's theorem, for every bounded continuous $\phi$,
    \begin{equation} \label{eq: duhamel weak}
        \langle u(\cdot,t),\phi\rangle=-\frac{\alpha}{\lambda}\int_0^t\qbheL(s)\bigl\langle \mu,P_{t-s}\phi\bigr\rangle\,ds,
    \end{equation}
    where $(P_\tau\phi)(\mathbf{y}):=\int\Phi(\mathbf{x}-\mathbf{y},\tau)\phi(\mathbf{x})\,d\mathbf{x}$ for $\tau>0$ and $P_0\phi:=\phi$. Again because $\Phi(\cdot,\tau)$ is a probability density, $\|P_\tau\phi\|_\infty\le\|\phi\|_\infty$, and $P_\tau\phi\to\phi$ pointwise as $\tau\to0$. Now let $\phi\in C_c^\infty(\mathbb{R}^3)$. Let $\Phi_0(\mathbf{z},\tau):=(4\pi\alpha\tau)^{-3/2}e^{-\|\mathbf{z}\|^2/(4\alpha\tau)}$ be the heat kernel, which satisfies $\partial_\tau\Phi_0=\alpha\Delta\Phi_0$. Writing $(P_\tau\phi)(\mathbf{y})=\int\Phi_0(\mathbf{z},\tau)\phi(\mathbf{y}+\mathbf{v}\tau+\mathbf{z})\,d\mathbf{z}$ and integrating by parts gives $\partial_\tau P_\tau\phi=P_\tau(A\phi)$ for $\tau>0$ with $A\phi:=\mathbf{v}\cdot\nabla\phi+\alpha\Delta\phi\in C_c^\infty(\mathbb{R}^3)$. Hence $P_\tau\phi=\phi+\int_0^\tau P_\sigma(A\phi)\,d\sigma$ pointwise, with $|P_\sigma(A\phi)|\le\|A\phi\|_\infty$. Inserting this into \eqref{eq: duhamel weak}, substituting $r=s+\sigma$ in the inner integral, exchanging the order of integration (Fubini, since $|\qbheL(s)\langle\mu,P_{r-s}(A\phi)\rangle|\le\|A\phi\|_\infty\,2L\Nb\,|\qbheL(s)|$) and applying \eqref{eq: duhamel weak} at time $r$ to the test function $A\phi$ gives
    \begin{multline*}
        \langle u(\cdot,t),\phi\rangle=-\frac{\alpha}{\lambda}\int_0^t\qbheL(s)\langle \mu,\phi\rangle\,ds-\frac{\alpha}{\lambda}\int_0^t\!\!\int_0^r\qbheL(s)\bigl\langle \mu,P_{r-s}(A\phi)\bigr\rangle\,ds\,dr\\
        =-\frac{\alpha}{\lambda}\int_0^t\qbheL(s)\langle \mu,\phi\rangle\,ds+\int_0^t\langle u(\cdot,r),A\phi\rangle\,dr .
    \end{multline*}
    Both integrands are integrable, since $|\qbheL(s)\langle \mu,\phi\rangle|\le\|\phi\|_\infty\,2L\Nb\,|\qbheL(s)|$ and $|\langle u(\cdot,r),A\phi\rangle|\le\|A\phi\|_\infty\|u(\cdot,r)\|_{L^1}$. Hence $t\mapsto\langle u(\cdot,t),\phi\rangle$ is absolutely continuous with derivative $\frac{\alpha}{\lambda}\langle f(\cdot,t),\phi\rangle+\langle u(\cdot,t),A\phi\rangle$, where $\langle f(\cdot,t),\phi\rangle=-\qbheL(t)\langle\mu,\phi\rangle$, for a.e.\ $t\in(0,T)$. Multiplying by $\zeta\in C_c^\infty(0,T)$ and integrating by parts in $t$ gives \eqref{eq: weak formulation}. Finally, $u(\cdot,0)=0$, and for every $\phi\in\mathscr S(\mathbb{R}^3)$ the map $t\mapsto\langle u(\cdot,t),\phi\rangle$ is continuous by \eqref{eq: duhamel weak} and dominated convergence, so $u\in C([0,T];\mathscr S'(\mathbb{R}^3))$.
\end{proof}

\begin{lemma}[MFLS representation] \label{lem: mfls representation}
    Let $\qbheL\in L^1_{\mathrm{loc}}([0,\infty))$ and let $f$ be the source \eqref{eq: source f}. At every point $\mathbf{x}=(\mathbf{x}_{xy},z)$ with $\mathbf{x}_{xy}\neq\mathbf{p}^k$ for all $k=1,\dots,\Nb$, the function $u(\cdot,t)$ from \eqref{eq: duhamel} is given by \eqref{eq: mfls}, and the convolution integrals in \eqref{eq: mfls} converge absolutely.
\end{lemma}
\begin{proof}
    Since $v_z=0$, $\Phi$ factorizes as $\Phi(\mathbf{x},\tau)=\Phi_H(\mathbf{x}_{xy}-\mathbf{v}_{xy}\tau,\tau)\,\Phi_V(z,\tau)$ with $\mathbf{v}_{xy}=(v_x,v_y)$, $\Phi_H(\mathbf{d},\tau)=e^{-\|\mathbf{d}\|^2/(4\alpha\tau)}/(4\pi\alpha\tau)$ and $\Phi_V(z,\tau)=e^{-z^2/(4\alpha\tau)}/\sqrt{4\pi\alpha\tau}$.
    For the line source $\mu_k$, the convolution is therefore the product
    \[
        (\Phi(\cdot,\tau)*\mu_k)(\mathbf{x}) = \Phi_H(\mathbf{x}_{xy}-\mathbf{p}^k-\mathbf{v}_{xy}\tau,\tau)\int_{\mathbb{R}}\Phi_V(z-z',\tau)\big(\mathbb{1}_{[0,L]}-\mathbb{1}_{[-L,0]}\big)(z')\,dz' .
    \]
    With $\operatorname{erf}z=\frac{2}{\sqrt\pi}\int_0^ze^{-t^2}\,dt$ \cite[Eq.~7.2.1]{DLMF}, the substitution $\xi=(z-z')/(2\sqrt{\alpha\tau})$ gives $\int_a^b\Phi_V(z-z',\tau)\,dz'=\frac12\big[\operatorname{erf}\big((z-a)/(2\sqrt{\alpha\tau})\big)-\operatorname{erf}\big((z-b)/(2\sqrt{\alpha\tau})\big)\big]$, so that the $z'$-integral above equals
    \[
        \tfrac12\Big[\operatorname{erf}\tfrac{z}{2\sqrt{\alpha\tau}}+\operatorname{erf}\tfrac{L-z}{2\sqrt{\alpha\tau}}\Big]-\tfrac12\Big[\operatorname{erf}\tfrac{L+z}{2\sqrt{\alpha\tau}}-\operatorname{erf}\tfrac{z}{2\sqrt{\alpha\tau}}\Big]=\mathcal{Z}(z,\tau,L).
    \]
    Let $\mathbf{d}=\mathbf{x}_{xy}-\mathbf{p}^k$ and $r=\|\mathbf{d}\|$. From $\|\mathbf{d}-\mathbf{v}_{xy}\tau\|^2=r^2-2\tau\langle\mathbf{d},\mathbf{v}_{xy}\rangle+\|\mathbf{v}\|^2\tau^2$,
    \[
        \frac{\alpha}{\lambda}\Phi_H(\mathbf{d}-\mathbf{v}_{xy}\tau,\tau)=\frac{1}{4\pi\lambda}\cdot\frac1\tau e^{-\frac{r^2}{4\alpha\tau}}\cdot e^{\frac{\langle\mathbf{d},\mathbf{v}\rangle}{2\alpha}}\cdot e^{-\frac{\|\mathbf{v}\|^2\tau}{4\alpha}}=C_\lambda\,\mathcal{R}(r,\tau)\,\mathcal{S}_\mathbf{v}(\mathbf{x}_{xy},\mathbf{p}^k)\,\mathcal{W}_\mathbf{v}(\tau).
    \]
    Combining, we get $\frac{\alpha}{\lambda}(\Phi(\cdot,\tau)*\mu_k)(\mathbf{x})=\mathcal{S}_\mathbf{v}(\mathbf{x}_{xy},\mathbf{p}^k)\,\mathcal{G}(r,z,\tau,L)$. Inserting $f=-\qbheL\sum_k\mu_k$ into \eqref{eq: duhamel} and substituting $\tau=t-s$ gives \eqref{eq: mfls}. For $r>0$ the integrand is absolutely integrable, because $\mathcal{G}(r,z,\cdot,L)$ is bounded for every $z\in\mathbb{R}$ and $\qbheL\in L^1(0,t)$. Its factors are bounded: for $a,c>0$,
    \begin{equation} \label{eq: sup bound}
        \sup_{\tau>0}\tau^{-a}e^{-c/\tau}=\Bigl(\frac{a}{ec}\Bigr)^a,
    \end{equation}
    so $\mathcal{R}(r,\cdot)\le4\alpha/(er^2)$, further $0<\mathcal{W}_\mathbf{v}\le1$, and $|\mathcal{Z}(z,\tau,L)|\le\int_{-L}^{L}\Phi_V(z-z',\tau)\,dz'\le\int_{\mathbb{R}}\Phi_V(z-z',\tau)\,dz'=1$, since $|\mathbb{1}_{[0,L]}-\mathbb{1}_{[-L,0]}|=\mathbb{1}_{[-L,L]}$ almost everywhere.

    The restriction $r>0$ is essential: for $r=0$ we have $\mathcal{R}=1/\tau$ and $\mathcal{Z}\to\pm1$ as $\tau\downarrow0$ for $0<\pm z<L$, so the integrand in \eqref{eq: mfls} behaves like $1/\tau$ near $\tau=0$ (the logarithmic singularity of the line source). On compact subsets of $\mathbb{R}^3\setminus\bigcup_k(\{\mathbf{p}^k\}\times\mathbb{R})$, the bounds above are uniform in $\mathbf{x}$, so by dominated convergence \eqref{eq: mfls} is continuous in $\mathbf{x}$ there and gives the canonical representative of the $L^1$ class $u(\cdot,t)$.
\end{proof}

\begin{lemma}[Uniqueness] \label{lem: uniqueness}
    Problem \eqref{eq: advection diffusion} has at most one solution in the sense of Definition~\ref{def: solution}.
\end{lemma}
\begin{proof}
    Let $w$ be the difference of two solutions. Then $w\in C([0,T];\mathscr S'(\mathbb{R}^3))$, $w(0)=0$, and $w$ satisfies \eqref{eq: weak formulation} with $f=0$. Since $t\mapsto\langle w(t),\phi\rangle$ is continuous on $[0,T]$ for every $\phi\in\mathscr S(\mathbb{R}^3)$, the family $\{w(t)\}_{t\in[0,T]}$ is pointwise bounded on the Fr\'echet space $\mathscr S(\mathbb{R}^3)$. It is therefore equicontinuous \cite[Thm.~2.6]{Rudin1991FA}: there are $C\ge0$ and $m\in\mathbb{N}$ with
    \[
        |\langle w(t),\phi\rangle|\le C\sum_{|\boldsymbol{\mu}|\le m}\sup_{\mathbf{y}\in\mathbb{R}^3}(1+\|\mathbf{y}\|)^m\,|\partial^{\boldsymbol{\mu}}\phi(\mathbf{y})| \qquad\text{for all } t\in[0,T],\ \phi\in\mathscr S(\mathbb{R}^3).
    \]
    Fix $\omega\in C_c^\infty(\mathbb{R}^3)$ and set $w_\omega(\mathbf{x},t):=\langle w(t),\omega(\mathbf{x}-\cdot)\rangle$. Since $\mathbf{x}\mapsto\omega(\mathbf{x}-\cdot)$ is smooth with values in $\mathscr S(\mathbb{R}^3)$, the estimate above shows that $w_\omega$ is smooth in $\mathbf{x}$ with $\partial^{\boldsymbol{\mu}}_{\mathbf{x}}w_\omega(\mathbf{x},t)=\langle w(t),(\partial^{\boldsymbol{\mu}}\omega)(\mathbf{x}-\cdot)\rangle$ and that, for fixed $t$, all these derivatives are continuous in $\mathbf{x}$. Moreover, they are jointly continuous on $\mathbb{R}^3\times[0,T]$: let $(\mathbf{x}_n,t_n)\to(\mathbf{x},t)$, $\psi_n:=(\partial^{\boldsymbol{\mu}}\omega)(\mathbf{x}_n-\cdot)$ and $\psi:=(\partial^{\boldsymbol{\mu}}\omega)(\mathbf{x}-\cdot)$. Then $\psi_n\to\psi$ in $\mathscr S(\mathbb{R}^3)$, since all $\psi_n$ are supported in a fixed compact set and depend smoothly on $\mathbf{x}_n$. Hence $|\langle w(t_n),\psi_n\rangle-\langle w(t),\psi\rangle|\le|\langle w(t_n),\psi_n-\psi\rangle|+|\langle w(t_n)-w(t),\psi\rangle|$. The first term is bounded by $C\sum_{|\boldsymbol{\mu}'|\le m}\sup_{\mathbf{y}}(1+\|\mathbf{y}\|)^m|\partial^{\boldsymbol{\mu}'}(\psi_n-\psi)(\mathbf{y})|\to0$ by the estimate above, which is uniform in $t$. The second term tends to $0$ because $t\mapsto\langle w(t),\psi\rangle$ is continuous. With $1+\|\mathbf{y}\|\le(1+\|\mathbf{x}\|)(1+\|\mathbf{x}-\mathbf{y}\|)$, it also gives $|\partial^{\boldsymbol{\mu}}_{\mathbf{x}}w_\omega(\mathbf{x},t)|\le C_{\boldsymbol{\mu}}(1+\|\mathbf{x}\|)^m$ uniformly in $t\in[0,T]$. Testing \eqref{eq: weak formulation} with $\phi=\omega(\mathbf{x}-\cdot)\in C_c^\infty(\mathbb{R}^3)$, for which $\nabla\phi=-(\nabla\omega)(\mathbf{x}-\cdot)$ and $\Delta\phi=(\Delta\omega)(\mathbf{x}-\cdot)$, gives $\partial_t w_\omega(\mathbf{x},\cdot) = -\mathbf{v}\cdot\nabla_{\mathbf{x}}w_\omega(\mathbf{x},\cdot)+\alpha\Delta_{\mathbf{x}}w_\omega(\mathbf{x},\cdot)$ in $\mathcal D'(0,T)$ for every $\mathbf{x}\in\mathbb{R}^3$. The right-hand side is continuous on $\mathbb{R}^3\times[0,T]$, so $w_\omega(\mathbf{x},\cdot)\in C^1([0,T])$, the equation holds pointwise, and $w_\omega(\cdot,0)=0$. The function $W(\mathbf{y},s):=w_\omega(\mathbf{y}+\mathbf{v}s/\alpha,\,s/\alpha)$ therefore solves $\partial_sW=\Delta W$ in $\mathbb{R}^3\times(0,\alpha T]$ with $W(\cdot,0)=0$. It belongs to $C^2_1(\mathbb{R}^3\times(0,\alpha T])\cap C(\mathbb{R}^3\times[0,\alpha T])$ and grows at most polynomially in $\mathbf{y}$. By the uniqueness theorem for the Cauchy problem of the heat equation \cite[Sec.~2.3.3, Thm.~7]{Evans2010PDE}, $W\equiv0$, i.e.\ $w_\omega\equiv0$.

    Finally, for $\phi\in C_c^\infty(\mathbb{R}^3)$ choose $\omega:=\phi(-\,\cdot)$. Then $\langle w(t),\phi\rangle=w_\omega(\mathbf{0},t)=0$, and since $C_c^\infty(\mathbb{R}^3)$ is dense in $\mathscr S(\mathbb{R}^3)$, $w(t)=0$ for all $t\in[0,T]$.
\end{proof}

\begin{proof}[Proof of Theorem~\ref{thm: existence}]
    Since $\qbheL\in L^1(0,T)$, Lemma~\ref{lem: duhamel solves} shows that \eqref{eq: duhamel} defines a solution with $\|u(\cdot,t)\|_{L^1}\le\frac{\alpha}{\lambda}\,2L\Nb\int_0^t|\qbheL|$; the second expression in \eqref{eq: duhamel} follows from $f(\cdot,s)=-\qbheL(s)\,\mu$. By Lemma~\ref{lem: uniqueness}, this solution is the only one. Lemma~\ref{lem: mfls representation} gives the representation \eqref{eq: mfls}.
\end{proof}

%% file: sections/Appendix_closed_form_V.tex
This appendix evaluates the integrated vertical response $\mathcal{Y}$ from \eqref{def: mathcalV} in closed form. Throughout, $\tau>0$, $L>0$ and $\vartheta := L/(2\sqrt{\alpha\tau}) > 0$.
Recall $\mathcal{Y}(\tau;L)=\int_0^L\mathcal{Z}(z,\tau,L)\,w_L(z)\,dz$ with $w_L(z)=N_W(L-z)/D_W(L)$, where $N_W(z)=C_1e^{M_1z}+C_2e^{M_2z}$ and $D_W$ are given in \eqref{sys: constants}, and $\mathcal{Z}$ in \eqref{eq: kernel Z}. Since $\gamma^2-\beta^2=\beta_1\beta_2+\beta_{12}(\beta_1+\beta_2)>0$, the exponents $M_{1,2}=\beta\pm\gamma$ satisfy $M_1>0>M_2$. In particular $M_1,M_2\neq0$.

\begin{proposition}[Closed form of $\mathcal{Y}$] \label{prop: closed V}
For $\vartheta>0$ and $\sigma\in\mathbb{R}\setminus\{0\}$ let
\begin{equation*}
    \mathcal{I}(\vartheta,\sigma) := \frac{1}{2}\int_0^1 K_0(x,\vartheta)\,e^{\sigma(1-x)}\,dx, \qquad K_0(x,\vartheta) := \operatorname{erf}(\vartheta(1-x)) + 2\operatorname{erf}(\vartheta x) - \operatorname{erf}(\vartheta(1+x)).
\end{equation*}
Then
\begin{equation} \label{eq: V closed}
    \mathcal{Y}(\tau;L) = \frac{L}{D_W(L)}\Bigl[C_1\,\mathcal{I}(\vartheta,M_1L)+C_2\,\mathcal{I}(\vartheta,M_2L)\Bigr]
\end{equation}
and
\begin{equation} \label{eq: I closed}
    \mathcal{I}(\vartheta,\sigma) = \frac{1}{2\sigma}\Bigl[\operatorname{erf}(2\vartheta)-2\operatorname{erf}(\vartheta)+\mathcal{Q}(0,\vartheta,\sigma)-2\mathcal{Q}(1,\vartheta,\sigma)+\mathcal{Q}(2,\vartheta,\sigma)\Bigr],
\end{equation}
where
\begin{equation*}
    \mathcal{Q}(n,\vartheta,\sigma) := \exp\!\Bigl(\frac{\sigma^2}{4\vartheta^2}+n\sigma\Bigr)\left[\operatorname{erf}\Bigl(\frac{\sigma}{2\vartheta}+(n-1)\vartheta\Bigr)-\operatorname{erf}\Bigl(\frac{\sigma}{2\vartheta}+n\vartheta\Bigr)\right].
\end{equation*}
\end{proposition}

\begin{proof}
\emph{Step 1: Reduction to $\mathcal{I}$.} With $z=Lx$ we have $\mathcal{Z}(Lx,\tau,L)=\frac12K_0(x,\vartheta)$ and $N_W(L-Lx)=\sum_{l=1}^2C_le^{M_lL(1-x)}$, hence $D_W(L)\,\mathcal{Y}(\tau;L) = L\sum_{l=1}^2C_l\,\mathcal{I}(\vartheta,M_lL)$, which is \eqref{eq: V closed}.

\emph{Step 2: Evaluation of $\mathcal{I}$.} Split $2\mathcal{I}=J_1+2J_2-J_3$ with
\begin{equation*}
    J_1=\int_0^1\operatorname{erf}(\vartheta y)\,e^{\sigma y}dy,\qquad J_2=\int_0^1\operatorname{erf}(\vartheta x)\,e^{\sigma(1-x)}dx,\qquad J_3=\int_0^1\operatorname{erf}(\vartheta(1+x))\,e^{\sigma(1-x)}dx,
\end{equation*}
where we substituted $y=1-x$ in the first integral. The antiderivative of $\operatorname{erf}(az+b)\,e^{\delta z+c}$ \cite[1.5.1]{korotkov_integrals} gives
\begin{align*}
    J_1 &= \frac1\sigma\Bigl[e^{\sigma}\operatorname{erf}(\vartheta)+\mathcal{Q}(0,\vartheta,\sigma)\Bigr],\qquad
    J_2 = -\frac1\sigma\Bigl[\operatorname{erf}(\vartheta)+\mathcal{Q}(1,\vartheta,\sigma)\Bigr],\\
    J_3 &= -\frac1\sigma\Bigl[\operatorname{erf}(2\vartheta)-e^{\sigma}\operatorname{erf}(\vartheta)+\mathcal{Q}(2,\vartheta,\sigma)\Bigr].
\end{align*}

\emph{Step 3: Assembly.} Inserting these expressions into $2\mathcal{I}=J_1+2J_2-J_3$ cancels the terms with $e^{\sigma}\operatorname{erf}(\vartheta)$ and gives \eqref{eq: I closed}.
\end{proof}

\begin{remark}[Symmetric BHEs] \label{rem: symmetric}
For symmetric BHEs with $\Rdelta1=\Rdelta2$, such as 1U and 2U configurations, $\beta_1=\beta_2$ and hence $\beta=0$. Then $M_{1,2}=\pm\gamma$, $C_1=C_2=\nu\gamma$ and $D_W(L)=\cosh(\gamma L)+\nu\sinh(\gamma L)$. With $\sigma:=\gamma L$, \eqref{eq: V closed} becomes
\begin{equation} \label{eq: V symmetric}
    \mathcal{Y}(\tau;L)=\frac{\nu}{2\bigl(\cosh\sigma+\nu\sinh\sigma\bigr)}\sum_{n=0}^{2}\binom{2}{n}(-1)^n\Bigl[\mathcal{Q}(n,\vartheta,\sigma)-\mathcal{Q}(n,\vartheta,-\sigma)\Bigr].
\end{equation}
\end{remark}

%% file: sections/Appendix_closed_form_centroid.tex
This appendix proves Lemma~\ref{lemma: centroid}. We use the notation of the lemma and write, with the vertex indices understood cyclically,
\begin{equation*}
    \Delta\mathbf{P}_i := \mathbf{P}_{i+1}-\mathbf{P}_i, \qquad \ell_i := \|\Delta\mathbf{P}_i\|, \qquad \mathbf{n}_i := \frac{\Delta\mathbf{P}_i^\perp}{\ell_i} .
\end{equation*}
For a counter-clockwise oriented polygon, $\mathbf{n}_i$ is the outward unit normal of the edge $S_i$. By the definitions of $\mathbf{n}_i$ and $W_i$,
\begin{equation} \label{eq: edge identities}
    (\mathbf{a}\cdot\mathbf{n}_i)\,\ell_i = W_i, \qquad \mathbf{n}_i\,\ell_i = \Delta\mathbf{P}_i^\perp .
\end{equation}
Parametrizing the edge by $\mathbf{x}(t)=\mathbf{P}_i+t\,\Delta\mathbf{P}_i$, $t\in[0,1]$, we have $ds=\ell_i\,dt$ and $\mathbf{a}\cdot\mathbf{x}(t)=f_i+t\,\Delta f_i$. We define the edge integrals
\begin{equation} \label{eq: edge integrals}
    \mathcal{E}_{i,0} := \int_0^1 e^{f_i+t\Delta f_i}\,dt, \qquad \mathcal{E}_{i,1} := \int_0^1 t\,e^{f_i+t\Delta f_i}\,dt .
\end{equation}

\begin{proof}[Proof of Lemma~\ref{lemma: centroid}]
The weighted centroid is $\mathbf{g}=\mathbf{M}_V/M_V$ with the mass $M_V$ and the first moment $\mathbf{M}_V$,
\begin{equation*}
    M_V := \iint_V e^{\mathbf{a}\cdot\mathbf{x}}\,dA, \qquad \mathbf{M}_V := \iint_V \mathbf{x}\,e^{\mathbf{a}\cdot\mathbf{x}}\,dA .
\end{equation*}
We use the divergence theorem $\iint_V\nabla\cdot\mathbf{G}\,dA=\oint_{\partial V}\mathbf{G}\cdot\mathbf{n}\,ds$ \cite[Thm.~5.8.1]{evans2015measure}, where $\mathbf{n}=\mathbf{n}_i$ is constant on every edge $S_i$.

\emph{Step 1: The edge integrals.} With $\int_0^1e^{zt}\,dt=\phi_1(z)$ and, by integration by parts, $\int_0^1te^{zt}\,dt=\phi_1(z)-\phi_2(z)$ for $z=\Delta f_i$, the integrals \eqref{eq: edge integrals} are those of Lemma~\ref{lemma: centroid}.

\emph{Step 2: The mass.} The vector field $\mathbf{F}(\mathbf{x}):=\|\mathbf{a}\|^{-2}e^{\mathbf{a}\cdot\mathbf{x}}\,\mathbf{a}$ satisfies $\nabla\cdot\mathbf{F}=\|\mathbf{a}\|^{-2}(\mathbf{a}\cdot\mathbf{a})\,e^{\mathbf{a}\cdot\mathbf{x}}=e^{\mathbf{a}\cdot\mathbf{x}}$. The divergence theorem, the parametrization of the edges and \eqref{eq: edge identities} give
\begin{equation} \label{eq: mass}
    M_V=\sum_i\frac{\mathbf{a}\cdot\mathbf{n}_i}{\|\mathbf{a}\|^2}\int_{S_i}e^{\mathbf{a}\cdot\mathbf{x}}\,ds=\frac{1}{\|\mathbf{a}\|^{2}}\sum_i(\mathbf{a}\cdot\mathbf{n}_i)\,\ell_i\,\mathcal{E}_{i,0}=\frac{1}{\|\mathbf{a}\|^{2}}\sum_iW_i\,\mathcal{E}_{i,0}.
\end{equation}

\emph{Step 3: The first moment.} For $j\in\{1,2\}$ let $\mathbf{Q}_j(\mathbf{x}):=\|\mathbf{a}\|^{-2}e^{\mathbf{a}\cdot\mathbf{x}}\,(x_j\mathbf{a}-\mathbf{e}_j)$. By the product rule,
\begin{equation*}
    \nabla\cdot\mathbf{Q}_j=\frac{e^{\mathbf{a}\cdot\mathbf{x}}}{\|\mathbf{a}\|^{2}}\Bigl(\mathbf{a}\cdot(x_j\mathbf{a}-\mathbf{e}_j)+\nabla\cdot(x_j\mathbf{a})\Bigr)=\frac{e^{\mathbf{a}\cdot\mathbf{x}}}{\|\mathbf{a}\|^{2}}\bigl(x_j\|\mathbf{a}\|^2-a_j+a_j\bigr)=x_j\,e^{\mathbf{a}\cdot\mathbf{x}} .
\end{equation*}
The divergence theorem, applied to both components, yields
\begin{equation*}
    \mathbf{M}_V=\frac{1}{\|\mathbf{a}\|^2}\sum_i\int_{S_i}e^{\mathbf{a}\cdot\mathbf{x}}\bigl[(\mathbf{a}\cdot\mathbf{n}_i)\,\mathbf{x}-\mathbf{n}_i\bigr]\,ds=\frac{1}{\|\mathbf{a}\|^2}\sum_i\ell_i\int_0^1e^{f_i+t\Delta f_i}\bigl[(\mathbf{a}\cdot\mathbf{n}_i)(\mathbf{P}_i+t\,\Delta\mathbf{P}_i)-\mathbf{n}_i\bigr]\,dt .
\end{equation*}
Evaluating the integral with \eqref{eq: edge integrals} and using \eqref{eq: edge identities}, we obtain
\begin{equation*}
    \mathbf{M}_V=\frac{1}{\|\mathbf{a}\|^2}\sum_i\Bigl[W_i\bigl(\mathbf{P}_i\,\mathcal{E}_{i,0}+\Delta\mathbf{P}_i\,\mathcal{E}_{i,1}\bigr)-\Delta\mathbf{P}_i^\perp\,\mathcal{E}_{i,0}\Bigr]=\frac{1}{\|\mathbf{a}\|^2}\sum_i\bigl[W_i\,\mathbf{m}_i-\Delta\mathbf{P}_i^\perp\,\mathcal{E}_{i,0}\bigr],
\end{equation*}
because $\mathbf{P}_i\mathcal{E}_{i,0}+(\mathbf{P}_{i+1}-\mathbf{P}_i)\mathcal{E}_{i,1}=\mathbf{P}_i(\mathcal{E}_{i,0}-\mathcal{E}_{i,1})+\mathbf{P}_{i+1}\mathcal{E}_{i,1}=\mathbf{m}_i$.

\emph{Step 4: The centroid.} Since $M_V>0$, also $\sum_iW_i\mathcal{E}_{i,0}=\|\mathbf{a}\|^2M_V>0$ by \eqref{eq: mass}, and in $\mathbf{g}=\mathbf{M}_V/M_V$ the factor $\|\mathbf{a}\|^{-2}$ cancels, which gives the assertion.
\end{proof}

%% file: sections/Appendix_derivatives.tex
This appendix proves Propositions~\ref{prop:soil jacobian}--\ref{prop:fluid jacobian}.

\subsection{Proof of Proposition~\ref{prop:soil jacobian}} \label{subsec: proof soil}
Let $t\in(0,T]$. Since $\qbheL=\qbhe/L$, the soil temperature \eqref{eq: mfls} at $\mathbf{x}=(\mathbf{x}_{xy},L/2)$ is a parametric convolution of the \emph{$\boldsymbol{\theta}$-independent} load $\qbhe$,
\begin{equation*}
u(t;\boldsymbol{\theta})=\bigl(\qbhe\ast K_x\bigr)(t),\qquad K_x(\tau;\boldsymbol{\theta}):=-\frac1L\sum_{k=1}^{\Nb}\mathcal{S}_k\,\mathcal{G}_k(\tau),
\end{equation*}
so all dependence on $L$, including that of the load per unit length, lies in the kernel, where the chain rule applies.

\emph{Regularity.} On $\Theta_x$ we have $r_k=\|\mathbf{x}_{xy}-\mathbf{p}^k\|>0$, so $\boldsymbol{\theta}\mapsto r_k$ is smooth; $\mathcal{S}_k$, $\mathcal{R}$, $\mathcal{W}_{\mathbf{v}}$, $\Zmid$ and $1/L$ are smooth functions of their arguments for $\tau>0$. Hence all $\boldsymbol{\theta}$-derivatives of $K_x$ exist and are continuous on $(0,T]\times\Theta_x$.

\emph{Derivatives.} With $\mathcal{S}_k=\exp(\langle\mathbf{x}_{xy}-\mathbf{p}^k,\mathbf{v}\rangle/(2\alpha))$ we have $\nabla_{\mathbf{p}^k}\mathcal{S}_k=-\frac{\mathbf{v}_{xy}}{2\alpha}\mathcal{S}_k$, and $\nabla_{\mathbf{p}^k}r_k=-(\mathbf{x}_{xy}-\mathbf{p}^k)/r_k$. From $\partial_r\mathcal{R}=-\frac{r}{2\alpha\tau}\mathcal{R}$ we get by the chain rule
\begin{equation*}
\nabla_{\mathbf{p}^k}\mathcal{G}_k=\partial_r\mathcal{G}_{\text{mid}}\,\nabla_{\mathbf{p}^k}r_k=\frac{\mathbf{x}_{xy}-\mathbf{p}^k}{2\alpha\tau}\,\mathcal{G}_k,\qquad\text{hence}\qquad
\nabla_{\mathbf{p}^k}K_x=-\frac{\mathcal{S}_k}{L}\,\mathbf{b}_k\,\mathcal{G}_k .
\end{equation*}
Among the factors of $\mathcal{G}_k$, only $\Zmid(\tau,L)=\mathcal{Z}(L/2,\tau,L)$ depends on $L$, through both of its arguments. With $\sigma=2\sqrt{\alpha\tau}$, $\partial_z\mathcal{Z}=\frac{1}{\sqrt\pi\sigma}\bigl[-e^{-(L-z)^2/\sigma^2}+2e^{-z^2/\sigma^2}-e^{-(L+z)^2/\sigma^2}\bigr]$ and $\partial_L\mathcal{Z}=\frac{1}{\sqrt\pi\sigma}\bigl[e^{-(L-z)^2/\sigma^2}-e^{-(L+z)^2/\sigma^2}\bigr]$, so the chain rule at $z=L/2$ gives
\begin{equation*}
\partial_L\Zmid=\tfrac12\partial_z\mathcal{Z}+\partial_L\mathcal{Z}=\frac{3}{2\sqrt\pi\sigma}\Bigl[e^{-\frac{L^2}{4\sigma^2}}-e^{-\frac{9L^2}{4\sigma^2}}\Bigr]=\frac{3}{4\sqrt{\pi\alpha\tau}}\Bigl[e^{-\frac{L^2}{16\alpha\tau}}-e^{-\frac{9L^2}{16\alpha\tau}}\Bigr].
\end{equation*}
By the product rule, $\partial_LK_x=-\sum_k\mathcal{S}_k\bigl[\tfrac1L\partial_L\mathcal{G}_{\text{mid}}(r_k,\tau,L)-\tfrac1{L^2}\mathcal{G}_k\bigr]$ with $\partial_L\mathcal{G}_{\text{mid}}=\Lambda\,\partial_L\Zmid$.

\emph{Domination.} Let $\Theta_0\subset\Theta_x$ be compact. By continuity and compactness there are $0<r_0\le r_1$ and $0<L_1\le L_2$ with $r_k\in[r_0,r_1]$ for all $k$ and $L\in[L_1,L_2]$ on $\Theta_0$; moreover $\mathcal{S}_k\le e^{\|\mathbf{v}\|r_1/(2\alpha)}$ and $\|\mathbf{x}_{xy}-\mathbf{p}^k\|\le r_1$. By \eqref{eq: sup bound} with $c=r_0^2/(4\alpha)$, $\mathcal{R}\le4\alpha/(er_0^2)$ and $\tau^{-1}\mathcal{R}\le(8\alpha/(er_0^2))^2$, and with $s=\tfrac12$ and $c=L_1^2/(16\alpha)$, $\partial_L\Zmid$ is bounded. Further, $0<\mathcal{W}_{\mathbf{v}}\le1$ and $0\le\Zmid\le1$, since $\Zmid=\tfrac12(3\operatorname{erf}x-\operatorname{erf}3x)$ with $x=L/(4\sqrt{\alpha\tau})$ increases from $0$ to $1$ in $x$. Hence $\mathcal{G}_k$, $\tau^{-1}\mathcal{G}_k$ and $\partial_L\mathcal{G}_{\text{mid}}=\Lambda\,\partial_L\Zmid$ are bounded on $(0,T]\times\Theta_0$, hence so are $K_x$, $\partial_LK_x$ and $\nabla_{\mathbf{p}^k}K_x$. This is \eqref{eq:domination}.

Hence $u(t;\cdot)\in C^1(\Theta_x)$ and, by \eqref{eq: derivative convolution}, $D_{\boldsymbol{\theta}}u=\bigl(\qbhe\ast D_{\boldsymbol{\theta}}K_x\bigr)(t)$. Inserting the derivatives and using $\qbhe/L=\qbheL$ yields \eqref{eq: dp u} and
\begin{equation*}
\partial_Lu=-\sum_k\mathcal{S}_k\bigl(\qbheL\ast\partial_L\mathcal{G}_{\text{mid}}(r_k,\cdot,L)\bigr)(t)+\frac1L\sum_k\mathcal{S}_k\bigl(\qbheL\ast\mathcal{G}_k\bigr)(t) .
\end{equation*}
By \eqref{eq: mfls}, $\sum_k\mathcal{S}_k(\qbheL\ast\mathcal{G}_k)(t)=-u(t;\boldsymbol{\theta})$, so the last term equals $-u/L$, which is \eqref{eq: dL u}. $\qed$

\subsection{Proof of Proposition~\ref{prop:psi jacobian}} \label{subsec: proof psi}
As above, $\Psi_{\text{soil}}(t;\boldsymbol{\theta})=\bigl(\qbhe\ast K_\Psi\bigr)(t)$ with the $\boldsymbol{\theta}$-independent load $\qbhe$ and
\begin{equation*}
K_\Psi(\tau;\boldsymbol{\theta}):=\frac{\psi(L)}{L}\,\mathcal{H}_\Sigma(\tau;\boldsymbol{\theta}),
\end{equation*}
where $\mathcal{H}_\Sigma$ is defined in \eqref{eq: Psi_soil}.
The diagonal terms $k=\ell$ have $\mathcal{S}_{\mathbf{v}}=1$ and $r_{k,k}=r_b$, and depend on $L$ only; for $k\neq\ell$ we have $r_{k,\ell}>0$ on $\Theta_\Psi$, so $K_\Psi$ is smooth in $\boldsymbol{\theta}$.

\emph{Derivatives with respect to $L$.} Only $\psi(L)/L$ and $\mathcal{H}(r,\tau,L)=\Lambda(r,\tau)\mathcal{Y}(\tau;L)$ depend on $L$; by the product rule
\begin{equation*}
\partial_LK_\Psi=\Bigl(\frac{\psi'(L)}{L}-\frac{\psi(L)}{L^2}\Bigr)\mathcal{H}_\Sigma+\frac{\psi(L)}{L}\,\partial_L\mathcal{H}_\Sigma,\qquad \partial_L\mathcal{H}_\Sigma=\frac1{\Nb}\sum_{k,\ell}\mathcal{S}_{\mathbf{v}}(\mathbf{p}^k,\mathbf{p}^\ell)\,\Lambda(r_{k,\ell},\tau)\,\partial_L\mathcal{Y}(\tau;L).
\end{equation*}
Since $\Psi_{\text{soil}}=\psi(L)(\qbheL\ast\mathcal{H}_\Sigma)$ and $\psi(L)=\frac12+\frac1{2\nu}\coth(\gamma L)>0$, the first term produces $\bigl(\psi'/\psi-1/L\bigr)\Psi_{\text{soil}}$, and the second term gives the last term of \eqref{eq: dL Psi}.

\emph{Derivatives with respect to $\mathbf{p}^m$.} The diagonal terms $k=\ell$ do not depend on $\mathbf{p}$. For $k\neq\ell$ the terms $(k,\ell)$ and $(\ell,k)$ combine to $\hat{\mathcal{S}}_{\mathbf{v}}(\mathbf{p}^\ell,\mathbf{p}^k)\mathcal{H}(r_{k,\ell},\tau,L)$ with the symmetrized shift kernel $\hat{\mathcal{S}}_{\mathbf{v}}$ of Section~\ref{sec: lumping split}. Only the pairs $\{m,\ell\}$, $\ell\ne m$, depend on $\mathbf{p}^m$. By the chain rule,
\begin{equation*}
\nabla_{\mathbf{p}^m}\hat{\mathcal{S}}_{\mathbf{v}}=\frac{\mathbf{v}_{xy}}{\alpha}\sinh\Bigl(\frac{\langle\mathbf{p}^m-\mathbf{p}^\ell,\mathbf{v}\rangle}{2\alpha}\Bigr),\qquad \nabla_{\mathbf{p}^m}r_{m,\ell}=\frac{\mathbf{p}^m-\mathbf{p}^\ell}{r_{m,\ell}},\qquad \partial_r\mathcal{H}=-\frac{r}{2\alpha\tau}\mathcal{H},
\end{equation*}
so that $\nabla_{\mathbf{p}^m}\bigl[\hat{\mathcal{S}}_{\mathbf{v}}\mathcal{H}(r_{m,\ell},\tau,L)\bigr]=\frac{\mathcal{H}}{\alpha}\bigl[\mathbf{v}_{xy}\sinh a_{m\ell}-\cosh(a_{m\ell})\tfrac{\mathbf{p}^m-\mathbf{p}^\ell}{\tau}\bigr]$ (using $\hat{\mathcal{S}}_{\mathbf{v}}=2\cosh a_{m\ell}$); multiplication by $\psi(L)/(L\Nb)$ gives the kernel of \eqref{eq: dp Psi}.

\emph{Domination.} Let $\Theta_0\subset\Theta_\Psi$ be compact. Then $L\in[L_1,L_2]\subset(0,\infty)$ and $r_{k,\ell}\in[r_0,r_1]$ for $k\ne\ell$, where we may take $r_0\le r_b\le r_1$. The functions $\psi,\psi',\hat{\mathcal{S}}_{\mathbf{v}},\sinh(\cdot)$ are bounded on $\Theta_0$, and $\Lambda$, $\tau^{-1}\Lambda$ are bounded as in the proof of Proposition~\ref{prop:soil jacobian}. Since $D_W(L)=e^{\beta L}\bigl(\cosh(\gamma L)+\nu\sinh(\gamma L)\bigr)>0$ and $|\mathcal{Z}|\le1$ (proof of Lemma~\ref{lem: mfls representation}), $\mathcal{Y}$ is bounded. With $w_L=N_W(L-\cdot)/D_W(L)$, the Leibniz rule gives
\begin{equation} \label{eq: dV}
\begin{aligned}
\partial_L\mathcal{Y}(\tau;L)&=-\frac{D_W'(L)}{D_W(L)^2}\int_0^L\!\mathcal{Z}N_W(L-z)\,dz\\
&\quad+\frac1{D_W(L)}\Bigl[\mathcal{Z}(L,\tau,L)N_W(0)+\int_0^L\!\partial_L\mathcal{Z}\,N_W(L-z)\,dz+\int_0^L\!\mathcal{Z}\,N_W'(L-z)\,dz\Bigr],
\end{aligned}
\end{equation}
where $0\le\partial_L\mathcal{Z}$ (Section~\ref{subsec: proof soil}) is bounded by a Gaussian density in $z$, so $\int_0^L\partial_L\mathcal{Z}\,dz\le1$ uniformly in $\tau$. Hence $\partial_L\mathcal{Y}$ is bounded, and so are $\mathcal{H}$, $\tau^{-1}\mathcal{H}$ and $\partial_L\mathcal{H}=\Lambda\,\partial_L\mathcal{Y}$. Hence $K_\Psi$ and its first derivatives satisfy \eqref{eq:domination}. Hence $\Psi_{\text{soil}}(t;\cdot)\in C^1(\Theta_\Psi)$ with $D_{\boldsymbol{\theta}}\Psi_{\text{soil}}=\bigl(\qbhe\ast D_{\boldsymbol{\theta}}K_\Psi\bigr)(t)$, and rewriting $\qbhe/L=\qbheL$ gives \eqref{eq: dL Psi} and \eqref{eq: dp Psi}. $\qed$

\subsection{Proof of Proposition~\ref{prop:fluid jacobian}} \label{subsec: proof fluid}
By \eqref{eq:potential_decomposition}, $\Toutsys=\Psi_{\text{ground}}-\Psi_{\text{soil}}-\Psi_{\text{power}}$. The term $\Psi_{\text{soil}}$ is treated in Proposition~\ref{prop:psi jacobian}. The term $\Psi_{\text{power}}=(\psi(L)-1)\qbhe(t)/\Cf$ depends on $\boldsymbol{\theta}$ only through $\psi(L)$, since $\qbhe$ and $\Cf$ do not depend on $\boldsymbol{\theta}$; hence $\partial_L\Psi_{\text{power}}=\psi'(L)\qbhe(t)/\Cf$ and $\nabla_{\mathbf{p}}\Psi_{\text{power}}=0$. Finally, by \eqref{Psi: ground} and Section~\ref{subsec: ground}, $\Psi_{\text{ground}}(L)=\chi(L)\sum_{l=1}^2C_lU_l(L)$ with $U_l(L)=\int_0^LT_0(\xi)e^{M_l(L-\xi)}d\xi$, which does not depend on $\mathbf{p}$. The integrand is smooth in $(\xi,L)$, and the Leibniz rule gives $U_l'(L)=T_0(L)+M_lU_l(L)$. Together with $D_W>0$ (proof of Proposition~\ref{prop:psi jacobian}) and the product rule this proves the formula for $\Psi_{\text{ground}}'$. $\qed$

%% file: sections/nomenclature.tex
\section*{Nomenclature}
\setstretch{1.2}
\addcontentsline{toc}{section}{Nomenclature}

\begin{longtable}{p{2.2cm} p{10.0cm} p{2.8cm}}
\toprule
\textbf{Symbol} & \textbf{Description} & \textbf{Unit} \\
\midrule
\endhead

\bottomrule
\endfoot

\multicolumn{3}{l}{\textbf{Acronyms and Abbreviations}} \\
BHE & Borehole Heat Exchanger & -- \\
CVT & Centroidal Voronoi Tessellation & -- \\
EED & Earth Energy Designer & -- \\
FFT & Fast Fourier Transform & -- \\
GSHP & Ground Source Heat Pump & -- \\
MFLS & Moving Finite Line Source & -- \\[0.5em]

\multicolumn{3}{l}{\textbf{Latin Symbols}} \\
$\Cf$ & Fluid heat capacity rate of a single borehole ($\mdot\,\cp$) & \si{\watt\per\kelvin} \\
$C_\lambda$ & Prefactor $1/(4\pi\lambda)$ of the horizontal kernel $\Lambda$ & \si{\meter\kelvin\per\watt} \\
$c_p$ & Specific heat capacity of the heat carrier fluid & \si{\joule\per\kilogram\per\kelvin} \\
$D_W, N_W$ & Denominator and numerator of the weight $w_L$ & --, \si{\per\meter} \\
$f(\mathbf{x}, t)$ & Volumetric thermal source term & \si{\watt\per\cubic\meter} \\
$\mathcal{G}, \mathcal{G}_{\mathrm{mid}}$ & MFLS kernel $\Lambda\mathcal{Z}$ and its value at $z=L/2$ & \si{\kelvin\meter\per\watt\per\second} \\
$g_{\mathrm{geo}}$ & Geothermal temperature gradient & \si{\kelvin\per\meter} \\
$\mathbf{g}_k(\eta)$ & Generator/centroid of Voronoi cell $k$ for the sharpness $\eta$ & \si{\meter} \\
$\mathcal{H}$ & Depth-integrated soil kernel $\Lambda\mathcal{Y}$ & \si{\kelvin\meter\per\watt\per\second} \\
$\mathcal{H}_\Sigma$ & Pair-summed soil kernel in $\Psi_{\mathrm{soil}}$ & \si{\kelvin\meter\per\watt\per\second} \\
$L$ & Borehole length, equal for all BHEs & \si{\meter} \\
$L_{\mathrm{min}}, L_{\mathrm{max}}$ & Lower and upper bounds of the borehole length & \si{\meter} \\
$\dot{m}$ & Mass flow rate through a single borehole & \si{\kilogram\per\second} \\
$M_V, \mathbf{M}_V$ & Mass and first moment of a polygon $V$ with respect to the CVT density & --, \si{\meter} \\
$\mathbf{n}_{kj}, b_{kj}$ & Outward unit normal and offset of edge $j$ of cell $V_k$ & --, \si{\meter} \\
$\dot{m}_{\mathrm{tot}}$ & Total mass flow rate of the BHE field ($N_b\,\dot{m}$) & \si{\kilogram\per\second} \\
$N_b$ & Number of borehole heat exchangers & -- \\
$n_\tau, \tau_{\mathrm{lo}}, h$ & Number, first node and logarithmic spacing of the quadrature nodes for $\kappa_j$ & --, \si{\second}, -- \\
$\mathbf{p}^k$ & Planar position $(p_x^k, p_y^k)^\top$ of the $k$-th BHE & \si{\meter} \\
$q_{\mathrm{bhe}}$ & Heat extraction rate per borehole & \si{\watt} \\
$q_{\mathrm{bhe},L}$ & Heat extraction rate per unit borehole length & \si{\watt\per\meter} \\
$q_{\mathrm{sys}}$ & Heat extraction/rejection load of the BHE field & \si{\watt} \\
$q \ast K$ & Causal convolution $\int_0^t q(t-\tau)K(\tau)\,d\tau$ & -- \\
$R_b, R_a$ & Borehole and internal thermal resistance (multipole method) & \si{\meter\kelvin\per\watt} \\
$R_1^\Delta, R_2^\Delta, R_{12}^\Delta$ & Delta-circuit internal BHE thermal resistances & \si{\meter\kelvin\per\watt} \\
$\mathcal{R}$ & Radial diffusion kernel & \si{\per\second} \\
$r$ & Radial (horizontal) distance from a borehole axis & \si{\meter} \\
$r_b$ & Borehole radius & \si{\meter} \\
$r_{k,\ell}$ & Distance between BHE $k$ and BHE $\ell$ & \si{\meter} \\
$\mathcal{S}_{\mathbf{v}}, \hat{\mathcal{S}}_{\mathbf{v}}$ & Advection shift kernel and its symmetrized form & -- \\
$t$ & Time & \si{\second} \\
$t_{\mathrm{end}}$ & Simulated operation time & \si{\second} \\
$T_{\mathrm{surface}}$ & Undisturbed ground temperature at the surface & \si{\celsius} \\
$T_0(z)$ & Undisturbed ground temperature profile & \si{\celsius} \\
$\Ts(z,t)$ & Borehole wall temperature & \si{\celsius} \\
$\Ti(z,t)$,\newline $\To(z,t)$ & Inlet and outlet fluid temperature along the borehole & \si{\celsius} \\
$T_{\mathrm{in}}(t), T_{\mathrm{out}}(t)$ & Inlet and outlet temperature of a BHE & \si{\celsius} \\
$T_{\mathrm{out,sys}}$ & Mean outlet temperature of the BHE field & \si{\celsius} \\
$T_{\min}^{\text{fluid}}, T_{\max}^{\text{fluid}}$ & Lower and upper limit of the outlet temperature & \si{\celsius} \\
$\mathcal{T}^*_{\min}, \mathcal{T}^*_{\max}$ & Sets of critical times for the lower and upper outlet temperature limits & \si{\second} \\
$u(\mathbf{x}, t)$ & Temperature change in the subsurface (negative for net heat extraction) & \si{\kelvin} \\
$V_k(\eta), \hat{V}_k(\eta)$ & Voronoi cell of BHE $k$ and the shrunken cell for the sharpness $\eta$; $\hat V_k:=\hat V_k(\eta_*)$ in (P) & -- \\
$\mathbf{v}_{\mathrm D}$ & Darcy flux (groundwater) & \si{\meter\per\second} \\
$\mathbf{v}$ & Effective thermal advection velocity $(v_x, v_y, 0)^\top$ & \si{\meter\per\second} \\
$\mathcal{W}_\mathbf{v}$ & Advective decay factor & -- \\
$w_L$ & Weight $N_W(L-\cdot)/D_W(L)$ in the outlet temperature & \si{\per\meter} \\
$\mathbf{x}$ & Spatial coordinate $(x, y, z)^\top$ & \si{\meter} \\
$\mathbf{x}_{xy}$ & Planar spatial coordinate $(x, y)^\top$ & \si{\meter} \\
$\mathbf{x}_i, \mathcal{N}_{\mathrm{s}}$ & Soil sensors $(\mathbf{x}_i^{xy},L/2)$ and their index set & \si{\meter}, -- \\
$\mathcal{Y}$ & Integrated vertical response & -- \\
$z$ & Depth below the surface & \si{\meter} \\
$\mathcal{Z}, \mathcal{Z}_{\mathrm{mid}}$ & Vertical finite source kernel and its value at $z=L/2$ & -- \\[0.5em]

\multicolumn{3}{l}{\textbf{Greek Symbols}} \\
$\alpha$ & Soil thermal diffusivity & \si{\square\meter\per\second} \\
$\beta,\beta_j,\gamma,M_{1,2}$ & Coefficients of the local BHE model ($j\in\{1,2,12\}$) & \si{\per\meter} \\
$\chi$ & Ratio $\psi/D_W$ & -- \\
$\delta$ & Dirac delta distribution & -- \\
$\Delta_{\mathrm{bhe}}$ & Minimum distance between two boreholes & \si{\meter} \\
$\Delta_{\mathrm{neigh}}$ & Minimum distance of a borehole to the property boundary & \si{\meter} \\
$\Delta t$ & Time step & \si{\second} \\
$\Delta t_B$ & Lumping interval (block length) of $\Psi_{\mathrm{inter}}$ & \si{\second} \\
$\Delta T_{k,\ell}$ & Temperature change at BHE $k$ caused by BHE $\ell$ & \si{\kelvin} \\
$\Delta T_{\max}$ & Maximum permitted temperature change at the property boundary & \si{\kelvin} \\
$\varepsilon_{kj}$ & Margin by which edge $j$ of cell $V_k$ is shrunk ($\varepsilon=\Deltabhe/2$ on shared edges, $\Deltaneigh$ on $\partial\Omega$) & \si{\meter} \\
$\eta, \eta_*$ & Sharpness of the upstream bias of the CVT density and its accepted value & \si{\per\meter} \\
$\boldsymbol{\theta}$ & Design variables $(L,\mathbf{p}^1,\dots,\mathbf{p}^{N_b})$ & \si{\meter} \\
$\kappa_j$ & Weights of the discrete convolution & \si{\kelvin\meter\per\watt} \\
$\lambda$ & Soil thermal conductivity & \si{\watt\per\meter\per\kelvin} \\
$\Lambda$ & Horizontal kernel $C_\lambda\mathcal{R}\mathcal{W}_{\mathbf{v}}$ & \si{\kelvin\meter\per\watt\per\second} \\
$\nu$ & Coefficient $(\beta_1+\beta_2)/(2\gamma)$ of the local BHE model & -- \\
$(\rho c)_m, (\rho c)_w$ & Volumetric heat capacity of the saturated soil ($=\lambda/\alpha$) and of water & \si{\joule\per\cubic\meter\per\kelvin} \\
$\rho(\mathbf{x})$ & Density of the CVT (depends on $\eta$) & -- \\
$\tau$ & Time lag in the convolutions & \si{\second} \\
$\Phi$ & Fundamental solution of $\partial_t+\mathbf{v}\cdot\nabla-\alpha\Delta$ & \si{\per\cubic\meter} \\
$\phi_1, \phi_2$ & Normalized exponential functions & -- \\
$\psi$ & Outlet factor of the local BHE model, $\psi(z)=\frac{1}{2\nu}\coth(\gamma z)+\frac12$ & -- \\
$\Psi_{\mathrm{ground}}$ & Ground potential & \si{\kelvin} \\
$\Psi_{\mathrm{inter}}$ & Interaction potential between different BHEs & \si{\kelvin} \\
$\Psi_{\mathrm{power}}$ & Borehole resistance potential & \si{\kelvin} \\
$\Psi_{\mathrm{self}}$ & Self-interaction potential & \si{\kelvin} \\
$\Psi_{\mathrm{soil}}$ & Soil potential & \si{\kelvin} \\
$\Omega$ & Property & -- \\
$\Omega_{\mathrm{eff}}$ & Effective domain of admissible borehole positions in $\Omega$ & -- \\

\end{longtable}